\documentclass[11pt, english]{article}
\usepackage[T1]{fontenc}
\usepackage[latin9]{inputenc}
\usepackage{geometry}
\usepackage{color}
\usepackage{array}
\usepackage{float}
\usepackage{mathrsfs}
\usepackage{multirow}
\usepackage{amsmath}
\usepackage{amssymb}
\usepackage{graphicx}
\usepackage{rotfloat}
\usepackage{diagbox}
\usepackage{booktabs}
\usepackage{siunitx}
\usepackage{rotating}
\usepackage{makecell} 
\PassOptionsToPackage{normalem}{ulem}
\usepackage{ulem}
\usepackage{titlesec}
\usepackage{xurl}      
\usepackage{hyperref} 
\usepackage{longtable}
\usepackage{bm}
\usepackage{comment}
\usepackage{amsmath}        
\usepackage{amsfonts}       
\usepackage{algorithm}      
\usepackage{algorithmicx}   
\usepackage{algpseudocode} 
\usepackage{amsthm,enumitem}
\usepackage{graphicx}
\usepackage{subcaption}
\usepackage{float}
\newtheorem{prop}{Proposition}

\newtheorem{remark}{Remark}

\usepackage{setspace}
\usepackage{todonotes}
\usepackage{cite}
\usepackage{authblk}

\titleclass{\subsubsubsection}{straight}[\subsubsection]
\newcounter{subsubsubsection}[subsubsection]
\renewcommand\thesubsubsubsection{\thesubsubsection.\arabic{subsubsubsection}}
\titleformat{\subsubsubsection}
  {\normalfont\normalsize\itshape}{\thesubsubsubsection}{1em}{}
\titlespacing*{\subsubsubsection}{0pt}{\baselineskip}{0.5\baselineskip}

\title{\textbf{Leveraging Electric School Buses for Disaster Recovery: Optimizing Routing and Energy Scheduling via Branch-and-Price}}
\date{}
\author[1]{Sayed Hamid Hosseini Dolatabadi}
\author[2]{Yuchen Dong}
\author[1,*]{Tanveer Hossain Bhuiyan}
\author[2]{Bo Zeng}
\author[3]{Brian O'Neill}
\author[3]{Anthony Severson}

\affil[1]{Department of Mechanical, Aerospace, and Industrial Engineering, The University of Texas at San Antonio, TX, USA.}
\affil[2]{Department of Industrial Engineering, University of Pittsburgh, PA, USA.}
\affil[3]{San Antonio Fire Department, TX, USA.}

\affil[*]{Corresponding author: \href{mailto:tanveer.bhuiyan@utsa.edu}{tanveer.bhuiyan@utsa.edu};
Contributing authors: \href{mailto:sayedhamid.hosseinidolatabadi@utsa.edu}{sayedhamid.hosseinidolatabadi@utsa.edu}; \href{mailto:yud96@pitt.edu}{yud96@pitt.edu}; \href{mailto:bzeng@pitt.edu}{bzeng@pitt.edu}; \href{mailto:brian.o'neill@sanantonio.gov}{brian.o'neill@sanantonio.gov}; \href{mailto:anthony.severson@sanantonio.gov}{anthony.severson@sanantonio.gov}
}

\begin{document}

\maketitle

\begin{abstract}

Natural disasters threaten the resilience of power systems, causing widespread power outages that disrupt critical loads (e.g., hospitals) and endanger public safety. Compared to the conventional restoration methods that often have long response times, leveraging government-controlled electric school buses (ESBs) with large battery capacity and deployment readiness offers a promising solution for faster power restoration to critical loads during disasters while traditional maintenance is underway. Therefore, we study the problem of dispatching, routing, and scheduling a heterogeneous fleet of ESBs to satisfy the energy demand of critical isolated loads around disasters addressing the following practical aspects: combined transportation and energy scheduling of ESBs, multiple back-and-forth trips of ESBs between isolated loads and charging stations, and spatial-wise coupling among multiple ESB routes. We propose an efficient mixed-integer programming model for routing and scheduling ESBs, accounting for the practical aspects to minimize the total restoration cost over a planning horizon. We develop an efficient exact branch-and-price (B\&P) algorithm and a customized heuristic B\&P algorithm integrating dynamic programming and labeling algorithms. Numerical results based on a real case study of San Antonio disaster shelters and critical facilities demonstrate that our proposed exact B\&P and heuristic B\&P algorithms are computationally 121 and 335 times faster, respectively, than Gurobi. Using network sparsity to incorporate the limitation in shelter-ESB type compatibility in the model demonstrates that the total restoration cost increases, on average, by 207\% as the network becomes fully sparse compared to fully connected. The capacity utilization metric reflects that the proposed practical ESB routing and scheduling enables an ESB to meet the energy demand 4.5 times its effective usable capacity. Additionally, results demonstrate that the required fleet size increases by 225\% as the weather changes from normal to adverse.

\end{abstract}

\section{Introduction}\label{sec:Introduction}

The deployment of electric school buses (ESBs) has progressed significantly over the last ten years worldwide. In addition to offering clean and healthy transportation for students, they can also serve as mobile energy sources due to their large onboard batteries. Hence, they can be dispatched to key facilities to discharge electricity to (partially) meet energy needs. In this context, this paper addresses the operational challenge arising from utilizing them with vehicle-to-building (V2B) technology for fast power restoration of critical isolated loads during disaster scenarios. 
This paper presents a study that aims to find the optimal routing and scheduling for a fleet of heterogeneous ESBs to meet the energy demand of critical isolated loads (e.g., hospitals, emergency shelters, and senior homes) that are fully disconnected from the grid around disasters, minimizing both fixed and operational costs.

\subsection{Background and Motivation}\label{subsec:Motivation}

The increasing frequency and severity of extreme weather events and natural disasters, such as heatwaves, hurricanes, wildfires, and floods, pose significant risks to the stability and resiliency of power systems. The National Centers for Environmental Information of the National Oceanic and Atmospheric Administration reported that there were 20 climate-driven disasters in the U.S. in 2021, each resulting in economic losses exceeding \$1 billion~\cite{AEE2023}. Furthermore, according to the U.S. Federal Emergency Management Agency, more than 2,700 major disasters have been reported in the U.S. in the last two decades, with fires and severe storms being recognized as the most frequent incidents~\cite{FEMA2023}. These extreme weather events can lead to widespread power outages, affecting millions of people, and causing critical infrastructure, such as hospitals, shelters, and emergency response centers, to lose electricity~\cite{hutchinson20223}. The Department of Energy reported that severe weather events resulted in 679 power outages between 2003 and 2012. On average, each disaster affects at least 50,000 customers, costing the U.S. economy up to \$33 billion, annually~\cite{DOE2023}. 
As such, San Antonio, the second-largest city in Texas, also often experiences power outages due to multi-type extreme weather events; notably, the catastrophic winter storm in February 2021 resulted in nearly half of San Antonio without power and left about 400,000 homes and businesses without power for a week~\cite{TexasPubR2021}. The growing frequency and variety of large-scale public emergencies place an overwhelming burden on the San Antonio Fire Department (SAFD), the primary disaster response agency for the City of San Antonio.

Besides the substantial economic losses, power outages resulting from natural disasters pose serious health risks, especially for individuals with medical conditions who are disproportionately at high risk, such as asthma and Chronic Obstructive Pulmonary Disease, as they rely on electricity-dependent medical equipment. Therefore, the rapid restoration of critical facilities around disasters is an urgent need for sustaining life-saving essential services. Traditional restoration methods include dispatching repair crews to repair damaged infrastructure, and reconfiguring the power grid to reroute electricity around affected areas~\cite{shi2022co}. Despite being effective, these methods are often time-consuming, requiring extensive coordination and labor over an extended restoration period, particularly in severely affected regions. Clearly, such a long outage can impose severe health risks to vulnerable populations.

In response to the need for innovative alternatives, interest is rapidly growing in leveraging the available electric vehicles (EVs) as mobile power sources (MPSs) to support power system restoration via V2B technologies~\cite{advancingV2G2022, hutchinson20223}. However, light-duty EVs often face practical limitations when applied to large-scale disastrous situations due to their small battery capacities~\cite{hu2023resilience}. In contrast, medium-duty EVs (e.g., electric buses (EBs) and electric trucks) with their larger battery capacities (up to 660 kWh) can provide substantial power support to critical loads during extended outages, which has inspired researchers to integrate EBs for power restoration around disasters (e.g.,~\cite{wu2023distributed, gao2017resilience}).

However, personal EVs and EBs vary widely in availability, ownership, and readiness during crises, making it difficult to coordinate their use effectively during disasters. In contrast, ESBs, bolstered by their increased penetration, are government-owned assets with structured availability (e.g., schools are closed during disasters) for emergency response. The adoption of ESBs among school districts in the U.S. accelerated in recent years due to the significant federal and state government incentives, such as the Environmental Protection Agency (EPA) Clean School Bus Program~\cite{AEE2023, EPA2023}.
Similar to other states and cities in the U.S., the San Antonio Independent School District has received 18 ESBs to electrify its school bus fleet in early 2024~\cite{EnvironmentTexas2024}.

As publicly owned assets, ESBs can be readily coordinated and dispatched. Their robust design features---larger wheel size typically around 22.5"~\cite{a_z_bus_electric_school_buses}, larger seating capacity, and powerful battery---enable ESBs to operate effectively in adverse disaster conditions for both energy delivery and transportation. 

Nevertheless, despite the practical advantages of ESBs, little research is available on the efficient routing and scheduling of ESBs to restore power to critical isolated loads. While researchers studied the integration of EVs and EBs for disaster resiliency (e.g.,~\cite{erenouglu2022resiliency, li2021routing}), these investigations lack comprehensive and efficient routing and scheduling models considering heterogeneous EBs, route interdependency of EBs, multiple back-and-forth trips with partial energy delivery, and continuous SOC tracking. Addressing these practical issues or factors makes the operational problem computationally very difficult. Also, the ESB dispatching, scheduling, and routing solutions should be computed faster and with good accuracy to generate effective plans in practice. Yet, existing studies on EVs/EBs either simply used commercial solvers (e.g.,~\cite{su2022critical, alghamdi2023resilience}), or heuristic algorithms (e.g.,~\cite{xu2019enhancing,zhang2021coordinated}), or exact solution methods (e.g.,~\cite{lei2018routing, li2021preallocation}) that fail to provide faster ESB scheduling and routing solutions for practical applicability.
Therefore, the primary motivation of this study is to develop an efficient scheduling and routing model and fast algorithms to help emergency management agencies, such as SAFD, leverage ESBs for enhancing real-world disaster resilience.

\subsection{Related Literature}\label{subsec:Related_Literature}

As noted earlier, using EVs as MPSs to enhance distribution system resiliency has gained a nontrivial amount of attention recently. The existing studies fall into three key areas, including pre-positioning, dispatching, routing, and scheduling of (1) MPSs (e.g., mobile energy generators and mobile battery storage) and repair crews (RCs) using reconfiguration strategies~\cite{anokhin2021mobility, lei2019resilient}, (2) private and/or light-duty EVs~\cite{shin2016plug, froger2022electric}, and (3) medium-duty EBs~\cite{hu2023resilience, zhang2024strategic, fan2023optimization, ji2023optimal}.

MPSs in emergency power restoration have been extensively studied, often focusing on resource allocation, pre-positioning of MPSs, and coordination with RCs. Anokhin et al.~\cite{anokhin2021mobility} presented a mixed-integer linear programming (MILP) model for pre-positioning and dispatching different types of MPSs, including mobile energy storage systems and mobile energy generators (MEGs) in coordination with RCs for post-disaster restoration. Lei et al.~\cite{lei2019resilient} developed a non-convex mixed-integer nonlinear programming (MINLP) model for dispatching RCs and MPSs. Taheri et al.~\cite{taheri2020distribution} proposed a scenario-based stochastic MILP model for dispatching MEGs and RCs to enhance distribution system resiliency. Other studies have leveraged the bidirectional charging/discharging capability of privately owned or light-duty EVs and modeled strategic allocation and pre-positioning of EVs before disasters. Su et al.~\cite{su2022critical} proposed a mixed-integer second-order conic program model, solved by Gurobi, for EV-based power restoration after an accidental blackout. Mohammadi et al.~\cite{mohammadi2016healer} and Alghamdi et al.~\cite{alghamdi2023resilience} proposed nonlinear models for the optimal allocation and pre-positioning of EVs to enhance distribution system resiliency. Advanced solution methods, including multi-agent deep reinforcement learning \cite{wang2022multi}, column-and-constraint generation (C\&CG) \cite{cao2022resilience, lei2018routing}, and progressive hedging \cite{kim2018enhancing, shi2023optimal}, have also been employed to address these pre-positioning and scheduling problems for distribution system restoration. However, these studies primarily focused on pre-positioning or single-trip energy deliveries of MPSs and EVs. Some studies have also focused on dispatching EVs during restoration. Xu et al.~\cite{xu2019enhancing} introduced a nonlinear optimization framework, solved using a Genetic algorithm (GA), utilizing EVs to transport energy from unaffected areas to isolated households. By dispatching EVs with V2G capabilities, RCs, and renewable sources, Wang et al.~\cite{wang2024enhancing} proposed an MILP model to minimize the total downtime of affected communities and solved using Gurobi. However, they assumed that EVs follow a pre-defined schedule and did not consider energy constraints with back-and-forth trips. Ereno{\u{g}}lu et al.~\cite{erenouglu2022resiliency} presented a scenario-based stochastic MILP model, solved by CPLEX, to maximize the service time of critical loads by dispatching an on-call EV fleet. However, they focused on the coordination of EV fleets with the distribution system, rather than with isolated loads.

We note that the aforementioned studies either overlook the joint transportation and energy scheduling of EVs and MPSs or fail to capture the complexities of multiple back-and-forth trips of EVs. More importantly, the limited battery capacity and the uncertainty in availability and readiness of personal EVs during disasters can hinder their effective deployment in prolonged outage situations \cite{hu2023resilience, hutchinson20223}. In contrast, studies on medium-duty EBs highlight their practical advantages (i.e., larger battery capacity and longer range) for power system restoration \cite{AEE2023}. 
To enhance the distribution system resiliency using pre-hurricane allocation of EBs, studies proposed stochastic programming methods, including a stochastic MINLP model, solved by CPLEX 
\cite{gao2017resilience}, and a two-stage stochastic programming model, solved by C\&CG \cite{li2021preallocation}. Later, Hu et al.~\cite{hu2023resilience} and Zhang et al.~\cite{zhang2024strategic} presented nonlinear mathematical formulations for resiliency-oriented dispatching of EBs to minimize the load curtailment. However, the focus of these studies was mainly on the pre-disaster allocation or energy dispatching of EBs at fixed locations.

To improve the resiliency of the distribution system during disasters, Li et al.~\cite{li2021routing} developed an MILP model to optimize EB routing and scheduling alongside distribution system reconfiguration, solved by CPLEX. However, they assumed a pre-defined set of service and charging trips. Zhang et al.~\cite{zhang2021coordinated} developed a bi-level model solved using GA to optimize the network topology and traffic paths of EBs. However, they did not fully address the interdependent routing required to meet the demand of isolated loads. Li et al.~\cite{li2021resilient} proposed an MILP model (solved by CPLEX), where EBs, borrowed from a bus company, are dispatched to support islanded microgrids and minimize the EB rental cost. However, the authors only assumed the assignment of idle EBs to some charging stations (CSs) to provide power to the grid. To enhance system-wide restoration during extreme events, Wu et al.~\cite{wu2023distributed} presented a distributed restoration framework to address the need for privacy and strategic autonomy of EB companies. They proposed an MILP model that considers EB service requirements solved by Gurobi, and a bi-level framework for determining each EB's energy and service task, solved using an alternating direction multiplier method. However, they mostly focused on meeting the passenger demands of EBs and restoring the distribution system using the V2G capabilities of EBs at CSs, without delving into the details of the routing and scheduling of EBs to meet the demand of isolated loads.

While privately-owned EVs and EBs offer promising opportunities for power system restoration around disasters, their effectiveness is often hindered by their practical limitations (e.g., availability, ownership, and preparedness). In contrast, government-owned or -controlled ESBs with their organized management and accessibility, enable structured and efficient deployment during disasters \cite{AEE2023}. Nevertheless, despite the valuable insights from existing literature, several significant research gaps persist.
\textbf{First,} a noteworthy limitation is the absence of studies leveraging such ESBs with their practical advantages over private EVs and EBs that make ESBs a structured, effective, and reliable medium of power restoration to critical isolated loads. 
\textbf{Second,} existing studies did not adequately address the coordination of a fleet of heterogeneous MPSs (e.g., private EVs, EBs, MEGs, and government-owned ESBs). Using ESBs with distinct characteristics (cost, battery capacity, energy consumption rate, and technology compatibility) in a fleet requires an efficient approach to meet the varying energy demands of isolated loads. The limited availability of certain ESB types and compatibility constraints with specific isolated loads due to technology limitations add further complexity, an area that remains largely unexplored. \textbf{Third,} a critical gap is the lack of attention to the joint transportation and energy scheduling, route interdependency, and spatial-wise coupling among MPSs around disasters. Practical aspects of routing and scheduling MPSs, including coordinating multiple MPSs to simultaneously serve the same load, multiple back-and-forth trips for charging and discharging, continuous SOC tracking, and demand shifting over a planning horizon, have not been adequately addressed.
\textbf{Fourth,} existing studies typically rely on heuristic methods that lack an optimality guarantee, or off-the-shelf solvers (e.g., Gurobi) that suffer from computational difficulty in handling realistic-sized instances. Given the life-saving implications of ESBs for isolated load restoration around disasters, developing an exact and computationally efficient solution approach is an essential task to be addressed.

\subsection{Contributions}\label{subsec:Contribution}

To fill the aforementioned gaps, this paper makes the following key contributions. \textbf{First,} we develop a new and comprehensive MILP model to leverage  ESBs to achieve an organized, reliable, and fast power restoration. This model addresses the following practical and critical aspects in dispatching, routing, and scheduling an ESB fleet: (1) different ESB types, each with a distinct battery capacity, cost, energy consumption rate, and technology compatibility; (2) simultaneous consideration of both transportation and energy scheduling of ESBs; (3) comprehensive coordination among ESBs with multiple back-and-forth trips, introducing spatial-wise coupling (i.e., route interdependencies) among ESB routes; and (4) continuous SOC tracking of ESB batteries. \textbf{Second,} we analyze this model and derive a set of properties and valid inequalities for enhancements. \textbf{Third}, we develop an efficient branch-and-price (B\&P) algorithm (and its variants) to solve the MILP model, which significantly outperforms a commercial solver. \textbf{Fourth,} based on a real case study of disaster shelters and other critical facilities in San Antonio, TX, U.S., we generate a set of managerial insights for practitioners on (1) the performance of different solution approaches; (2) the effect of limited service capability of ESBs due to technology compatibility on the performance of different solution approaches, total cost, and required ESB fleet composition; (3) the effect of shelters' energy demand, shelters' service time, number of available ESBs, and disaster severity on the total cost and the required ESB fleet composition; and (4) the effect of the inconvenience fee of demand shifting on the percentage of shifted demand.

\section{Problem Description and Mathematical Formulation}\label{sec:Problem_Description_Math_Formulation}

\subsection{Problem Description}\label{sec:Problem_Description}

\begin{figure}[H]
    \centering
    \begin{subfigure}[b]{0.48\textwidth}
        \centering
        \includegraphics[width=\textwidth]{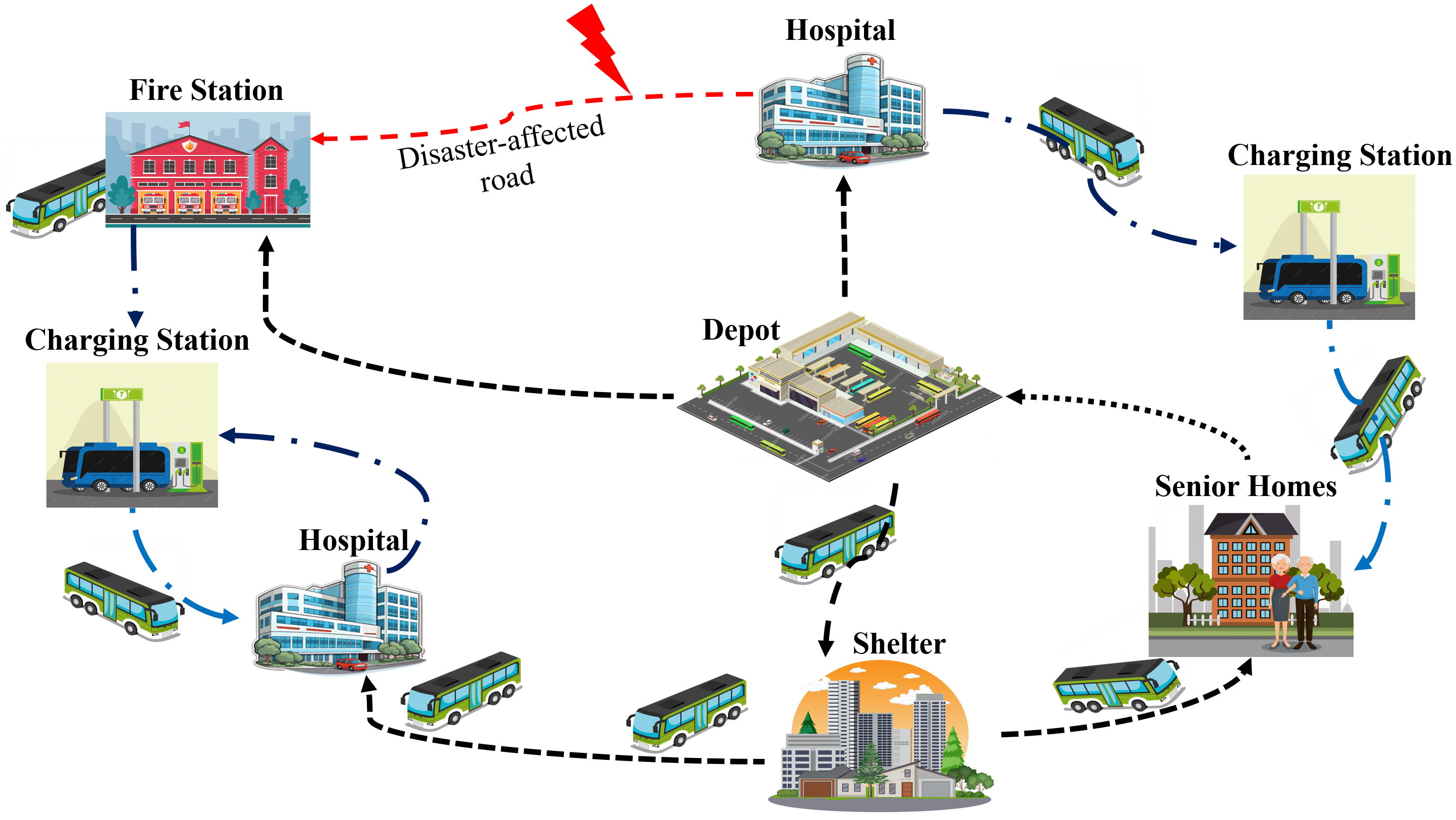}
        \caption{Routing and scheduling of ESBs.}
        \label{fig:Problem_Description_Figure}
    \end{subfigure}
    \begin{subfigure}[b]{0.48\textwidth}
        \centering
        \includegraphics[width=\textwidth]{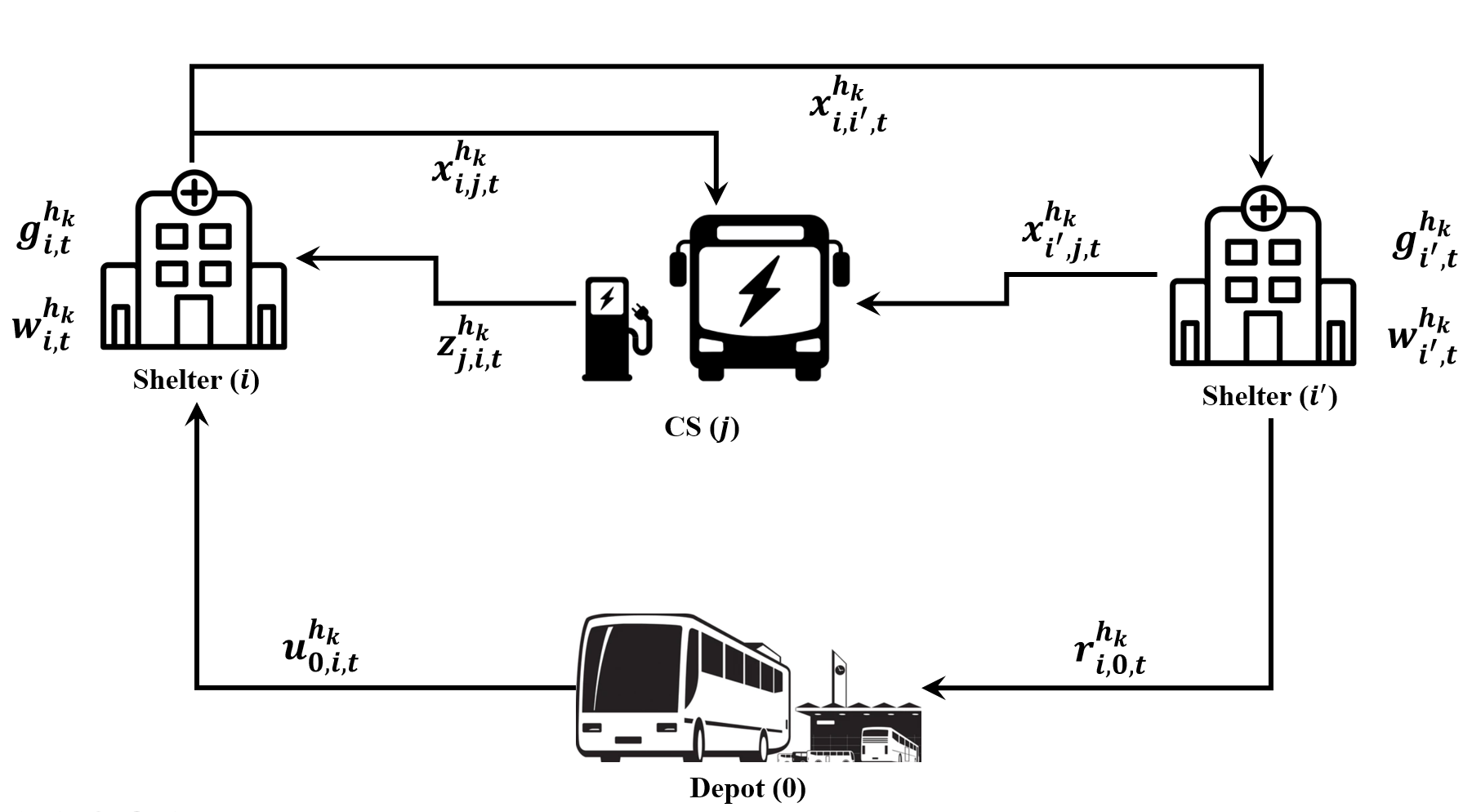}
        \caption{Model and decision variables.}
        \label{fig:Visualization-of-the-model}
    \end{subfigure}
    \caption{Visualization of the problem and model for routing and scheduling of ESBs.}
    \label{fig:Visulatization_of_problem_and_model}
\end{figure}

Consider employing ESBs to address power restoration for critical loads under a disaster, as shown in Figure \ref{fig:Problem_Description_Figure}, where the critical loads, such as hospitals, fire stations, and old homes (referred to as "shelters" throughout the remainder of the paper), are isolated from the distribution system due to a major grid outage contingency. Indeed, three types of nodes (i.e., locations) are involved in our network. (1) A \textit{central depot} node where a fleet of heterogeneous ESBs is stationed. Each ESB type has distinct characteristics, including cost, battery capacity, and energy consumption rate. (2) A set of \textit {isolated shelters} with fixed service times (i.e., the number of time slots it takes for an ESB to discharge to a shelter) that require power restoration. These shelters have time-varying energy demands that are assumed to be fixed and known \textit{a priori} over the planning horizon. However, some shelters may not accommodate certain ESB types to be plugged in due to technology, infrastructure, or spatial limitations, which restrict their service capability. (3) A set of \textit{available CSs} where ESBs can be charged for fixed service times. The CSs are assumed to have sufficient grid connectivity and capacity to meet the charging demand of all ESBs without limitation.

ESBs are dispatched from the depot to shelters as needed throughout the planning horizon, where they discharge their battery to meet the demand of shelters and must return to the depot by the end of the planning horizon. As the ESBs travel and serve shelters, their battery SOC levels deplete, requiring them to periodically visit CSs. Therefore, the ESBs make multiple back-and-forth trips between shelters and CSs to transport energy, ensuring that their battery SOC do not fall below the minimum acceptable level at any time. On the one hand, each ESB can serve multiple shelters before visiting a CS, depending on the energy demand of the shelters and the ESB's remaining SOC level. On the other hand, if the demand level of a shelter exceeds the battery capacity of an ESB, multiple ESBs are required to satisfy the demand of that shelter, introducing route interdependency and spatial-wise coupling among ESBs to ensure the energy demand of isolated shelters is satisfied.

The routing of ESBs accounts for the departure time and destination (i.e., when departing a source for a destination), and the availability status of each ESB at each shelter at each time slot (i.e., when an ESB is actively discharging to a shelter). As we have different types of nodes in the network, we use different routing actions that not only take into account the source-destination of the ESBs but also identify the timing of each trip (i.e., schedule). We consider a single arc between each pair of nodes in the network, where the distance and travel time on each arc are given \textit{a priori}. 
The energy schedule determines the amount of energy each ESB discharges to each shelter at each time slot during the planning horizon to satisfy the demand of all shelters while maintaining valid SOC levels of the ESBs at all times. All ESBs departing the depot or CSs have fully charged batteries, and any travel within the transportation network consumes energy. However, we assume that a fully-charged ESB (either departing from the depot or from a CS) can directly travel to a shelter without requiring a visit to an intermediate en-route CS.

Given all these aspects, the goal of the decision maker is to efficiently route and schedule the ESB fleet to minimize the total restoration cost while satisfying the energy demand of all shelters over the planning horizon. This problem inherently involves ESB fleet sizing (as the total cost is a function of the number of required ESBs) and represents an isolated bulk system restoration problem over a
planning horizon. The key outputs of the restoration program are the most efficient
routing plans, as well as charging, discharging, and dispatching schedules of ESBs that serve the energy demand of all shelters
at a minimum cost. We refer to this problem as an ESB routing and scheduling problem for isolated load restoration around disasters (\textit{ESBILR}).

\subsection{Mathematical Formulation}\label{sec:Mathematical_Formulation}

In this section, we present the mathematical formulation of the \textit{ESBILR}. We list the necessary sets, parameters, and variables supporting the formulation in Table \ref{tab:Sets}.

\begin{table}[h!]
\caption{\label{tab:Sets}Sets, parameters, and variables.}
\centering{}%
\begin{tabular}{c>{\raggedright}m{13cm}}
\hline 
$\textbf{Set}$ & $\textbf{Description}$\tabularnewline
\hline 
$\mathcal{I}^k$ & Set of shelters, indexed by $i$ that can be served by an ESB of type $k$\tabularnewline
$\mathcal{J}$ & Set of charging stations, indexed by $j$ \tabularnewline
$\mathcal{K}$ & Set of ESB types, indexed by $k$ \tabularnewline
$\mathcal{K}^i$ & Set of ESB types, indexed by $k$ that can serve shelter $i$ \tabularnewline
$\mathcal{H}^{k}$ & Set of available ESBs of type $k$ in the depot, indexed by $h_k$\tabularnewline
$\mathcal{T}$ & Planning horizon with time slots (discrete times $\Delta$), indexed
by $t$\tabularnewline
$\mathcal{Q}^k=\mathcal{I}^k\cup\mathcal{J}$ & Set of all CSs and shelters, indexed by $q$ that are accessible by ESBs of type $k$\tabularnewline
\hline 
$\textbf{Parameter}$ & $\textbf{Description}$\tabularnewline
\hline 
$C_{max}^{k}$ & Maximum charging capacity of ESBs of type $k$\tabularnewline
$C_{min}^{k}$ & Minimum required remaining energy in the battery of ESBs of type $k$ \tabularnewline
$P_{i,t}$ & Active power demand of shelter $i$ at time slot $t$\tabularnewline
$R_{a,b}^{k}$ & Energy consumption for ESBs of type $k$ to travel from node $a$ to node $b$
in the transportation network\tabularnewline
$T_{a,b}^{k}$ & Time required for ESBs of type $k$ to travel from node $a$ to node $b$ in
the transportation network\tabularnewline
$S_{j}$ & Service time of charging station $j$\tabularnewline
$S_{i}$ & Service time of shelter $i$\tabularnewline
$G_{min}^k$ & Minimum discharged energy from an ESB of type $k$ when actively serving (i.e., discharging energy) shelters\tabularnewline
$C_{inv}^k$ & Investment cost of using each ESB of type $k$ in the restoration program\tabularnewline
$C_{fx}$ & Fixed energy cost of ESBs (per unit electricity consumption cost)\tabularnewline
$H_{max}^{k}$ & Available number of ESBs of type $k$ at the depot\tabularnewline
$T_{1}$ & First time slot of the planning horizon\tabularnewline
$T_{2}$ & Last time slot of the planning horizon\tabularnewline
$F_i$ & Penalty cost for unsatisfied aggregated demand of shelter $i$\tabularnewline
\hline 
$\textbf{Variable}$ & $\textbf{Description}$\tabularnewline
\hline 
$x_{i,q,t}^{h_k}$ & Binary routing variable which is 1 if ESB $h$ of type $k$ is departing shelter $i$ for another shelter or
a charging station $q$ at time slot $t$, 0 otherwise\tabularnewline
$z_{j,i,t}^{h_k}$ & Binary routing variable which is 1 if ESB $h$ of type $k$ is departing charging station $j$ for shelter $i$
at time slot $t$, 0 otherwise\tabularnewline
$u_{0,i,t}^{h_k}$ & Binary routing variable which is 1 if ESB $h$ of type $k$ is departing the depot (indexed by 0) for shelter $i$ at time slot $t$, 0 otherwise\tabularnewline
$r_{i,0,t}^{h_k}$ & Binary routing variable which is 1 if ESB $h$ of type $k$ is departing shelter $i$ for the depot at time slot $t$, 0 otherwise\tabularnewline
$w_{i,t}^{h_k}$ & Binary routing variable which is 1 if ESB $h$ of type $k$ is serving shelter $i$ at time slot $t$, 0 otherwise\tabularnewline
$g_{i,t}^{h_k}$ & Continuous variable of the amount of energy from the battery of ESB $h$ of type $k$ discharged to shelter $i$ at time slot $t$\tabularnewline
$soc_{t}^{h_k}$ & Continuous variable of SOC of the battery of ESB $h$ of type $k$ at time slot $t$\tabularnewline
$l_i$ & Continuous variable for the unsatisfied aggregated demand of shelter $i$\tabularnewline
\hline
\end{tabular}
\end{table}

\vspace{-10pt}
\paragraph*{Description of Timing:}

As the planning horizon comprises discrete time slots, we define all routing and energy scheduling activities of ESBs in terms of time slots. Each time slot $t$ is defined with a specific time duration $\Delta$, and represents the beginning of the time interval $[t,\,t+\Delta)$. Each ESB can be in one of the following states at each time slot, including $(i)$ departing the depot for a shelter, $(ii)$ departing a shelter for another shelter or a CS, $(iii)$ departing a CS for a shelter, and $(iv)$ departing a shelter for the depot. 
An ESB starts its route by departing the depot, denoted by $u_{0,i,t}^{h_k}$, which represents that ESB $h$ of type $k$ departs the depot for shelter $i$ at time slot $t$ (i.e., the beginning of time interval $[t,\,t+\Delta)$). It means that before this time (i.e., $\forall t\in[0,\,t)\,or\,[0,\,t-\Delta]$) this ESB $h$ was staying at the depot. Then, ESB $h$ takes $T_{0,i}^k$ time slots to travel from the depot to shelter $i$. It means that at the beginning of time interval $[t+T_{0,i}^k,\,t+T_{0,i}^k+\Delta)$, equivalently at time slot $t+T_{0,i}^k$, the ESB $h$ reaches the shelter $i$. Therefore, the status of ESB $h$ changes to actively serving (i.e., discharging energy to) the shelter $i$, implying that $w_{i,t}^{h_k}=1$ from time slot $t+T_{0,i}^k$ until the beginning of time interval $[t+T_{0,i}^k+S_{i},\,t+T_{0,i}^k+S_{i}+\Delta)$, or equivalently from time slots $t+T_{0,i}^k$ to $t+T_{0,i}^k + S_i - 1$. Then, at the beginning of time interval $[t+T_{0,i}^{k}+S_{i},\,t+T_{0,i}^{k}+S_{i}+\Delta)$, ESB $h$ departs the shelter $i$ for another shelter or a CS. This means that at time slot $t+T_{0,i}^{k}+S_{i}$, $x_{i,q,t}^{h_k}=1$ for ESB $h$. Similarly, if an ESB $h$ departs a shelter $i$ at time slot $t$, i.e., $x_{i,q,t}^{h_k}=1$, traversing arc $i-q$ takes $T_{i,q}^{k}$ time slots for ESB $h$ to reach to either a CS or a shelter. Assuming the ESB $h$ arrives at a CS $j$ in this trip, this arrival happens at time slot $t+T_{i,j}^{k}$. Then, the ESB $h$ becomes fully charged within $S_{j}$ time slots, and at time slot $t+T_{i,j}^{k}+S_{j}$, it departs the CS $j$ for shelter $i$. Therefore, $z_{j,i,t}^{h_k}=1$ at time slot $t+T_{i,j}^{k}+S_{j}$. Similarly, if an ESB $h$ departs a CS $j$ for a shelter $i$ at time slot $t$ (i.e., $z_{j,i,t}^{h_k}=1$), this travel on arc $j-i$ takes $T_{j,i}^{k}$ time slots to reach to the shelter $i$. Therefore, at time slot $t+T_{j,i}^{k}$, the ESB $h$ arrives at shelter $i$ and discharges to this shelter $i$ for $S_{i}$ time slots (i.e., $w_{i,t}^{h_k}=1$ from time slots $t+T_{j,i}^{k}$ to $t+T_{j,i}^{k}+S_{i}-1$). Then, at time slot $t+T_{j,i}^{k}+S_{i}$, it departs the shelter $i$ for a CS or a shelter. Therefore, $x_{i,q,t}^{h_k}=1$ at time slot $t+T_{j,i}^{k}+S_{i}$.

Figure \ref{fig:Visualization-of-the-model} demonstrates a visualization of the model and variables associated with the nodes in the network. The detailed MILP formulation (referred to as \textit{ESB-BaM-MILP}) of the \textit{ESBILR} is presented in Eqs. (\ref{eq:Obj Function}) - (\ref{eq:Continuous_Variables}).

\paragraph*{Objective Function:}

\begin{align}
& min \qquad \Bigg\{\sum_{k\in\mathcal{K}}\sum_{h_k\in\mathcal{H}^k}\sum_{t\in\mathcal{T}}\sum_{i\in\mathcal{I}^k}\left(C_{inv}^{k}u_{0,i,t}^{h_k}\right)+ \sum_{i \in \mathcal{I}} F_i l_i +\label{eq:Obj Function}\\
&
C_{fx}\sum_{k\in\mathcal{K}}\sum_{h_k\in\mathcal{H}^k}\sum_{i\in\mathcal{I}^k}\sum_{t\in\mathcal{T}}\left(\sum_{q\in\mathcal{Q}^k}\left(x_{i,q,t}^{h_k}R_{i,q}^{k}\right) + \sum_{j\in\mathcal{J}}\left(z_{j,i,t}^{h_k}R_{j,i}^{k}\right) + \left(u_{0,i,t}^{h_k}R_{0,i}^{k}\right) + \left(r_{i,0,t}^{h_k}R_{i,0}^{k}\right)\right)\Bigg\}\nonumber
\end{align}

The objective function (\ref{eq:Obj Function}) seeks to minimize the investment cost of the required ESB fleet size (the first component), the penalty cost of unmet demand (the second component), and the total transportation energy consumption cost of ESBs (the third component).
As each ESB departs the depot at most once, the number of ESBs departing the depot represents the number of ESBs used (i.e., participated) in the critical load restoration program. Therefore, the first component of the objective function provides the minimum required ESB fleet size for the restoration program. If the available number of ESBs is not enough to satisfy the total energy demand of the shelters, a portion of the demand remains unsatisfied. Therefore, the second component computes the penalty cost of the total unmet demand.

\paragraph*{Critical Load Satisfaction Constraints:}

\begin{align}
&\sum_{k \in \mathcal{K}^i}\sum_{t\in\mathcal{T}}\sum_{h_k\in\mathcal{H}^k}g_{i,t}^{h_k} + l_i \geq \sum_{t\in\mathcal{T}}P_{i,t},\qquad\forall i\in\mathcal{I}\label{eq:Load Satisfaction - Weak Coupling}
\end{align}

Constraints (\ref{eq:Load Satisfaction - Weak Coupling}) ensure that the total energy demand of each shelter over the planning horizon (i.e., $\sum_{t\in\mathcal{T}}P_{i,t}$) is satisfied either solely by the aggregated discharged energy of ESBs, or along with incurring a penalty for any unmet demand (i.e., $l_i$) if the number of available ESBs is not enough to satisfy the demand.

\paragraph*{ESB Routing Constraints:}

\allowdisplaybreaks
\begin{align}
&\sum_{i\in\mathcal{I}^k}\sum_{t\in\mathcal{T}}u_{0,i,t}^{h_k}\leq 1,\forall k \in \mathcal{K}, h_k\in\mathcal{H}^k \label{eq:Routing - 1}\\
&\sum_{i\in\mathcal{I}^k}\sum_{q\in\mathcal{Q}^k}x_{i,q,t}^{h_k}+\sum_{i\in\mathcal{I}^k}w_{i,t}^{h_k}+\sum_{j\in\mathcal{\mathcal{J}}}\sum_{i\in\mathcal{I}^k}z_{j,i,t}^{h_k}+\sum_{i\in\mathcal{I}^k}r_{i,0,t}^{h_k}\leq\sum_{i\in\mathcal{I}^k}\sum_{\tau=T_1}^{t-1}u_{0,i,\tau}^{h_k},\forall k \in \mathcal{K}, h_k\in\mathcal{H}^k,\,t\in\mathcal{T} \label{eq:Routing - 2}\\
&\sum_{q\in\mathcal{Q}^k}x_{i,q,\left(t+T_{0,i}^{k}+S_{i}\right)}^{h_k}+r_{i,0,\left(t+T_{0,i}^{k}+S_{i}\right)}^{h_k}\geq u_{0,i,t}^{h_k}, \,\forall k \in \mathcal{K}, h_k\in\mathcal{H}^k,\,i\in\mathcal{I}^k,\,t+T_{0,i}^{k}+S_{i}<T_2\label{eq:Routing - 3}\\
&u_{0,i,t}^{h_k}=0, \,\forall k \in \mathcal{K}, h_k\in\mathcal{H}^k,\,i\in\mathcal{I}^k,\,t\in\mathcal{T},\,t+T_{0,i}^{k}+S_i\geq T_2\label{eq:Routing - 3-1}\\
&\sum_{i'\in\mathcal{I}^k}z_{j,i',\left(t+T_{i,j}^{k}+S_{j}\right)}^{h_k}\geq x_{i,j,t}^{h_k},\,\forall k \in \mathcal{K},h_k\in\mathcal{H}^k,\,i\in\mathcal{I}^k,\,j\in\mathcal{J},\,t\in\mathcal{T},\,t+T_{i,j}^{k}+S_{j}<T_2\label{eq:Routing 4 - 1}\\
& \sum_{q\in\mathcal{Q}^k}x_{i^{'},q,\left(t+T_{i,i^{'}}^{k}+S_{i^{'}}\right)}^{h_k}+r_{i',0,\left(t+T_{i,i^{'}}^{k}+S_{i^{'}}\right)}^{h_k}\geq x_{i,i^{'},t}^{h_k}, \label{eq:Routings 4 - 2}\\
& \,\forall k \in \mathcal{K}
,h_k\in\mathcal{H}^k,\,i,i^{'}\in\mathcal{I}^k\left(i'\neq i\right),\,t\in\mathcal{T},\,t+T_{i,i^{'}}^{k}+S_{i^{'}}<T_2\nonumber\\
& \sum_{q^{'}\in\mathcal{Q}^k}\sum_{\tau=t-S_{i}}^{t+T_{i,q}^{k}+S_{q}}\sum_{i^{'}\in\mathcal{I}^k\setminus\{i\}}x_{i^{'},q^{'},\tau}^{h_k} + \sum_{j^{'}\in\mathcal{J}}\sum_{\tau=t-S_{i}}^{t+T_{i,q}^{k}+S_{q}}\sum_{i^{'}\in\mathcal{I}^k}z_{j^{'},i^{'},\tau}^{h_k}\leq M_{1}\left(1-x_{i,q,t}^{h_k}\right),\label{eq:Routing - 5}\\
& \forall k \in \mathcal{K},
h_k\in\mathcal{H}^k,\,i\in\mathcal{I}^k,\, q\in\mathcal{Q}^k,\, t\in\mathcal{T},\,t+T_{i,q}^{k}+S_{q}<T_2\nonumber\\
&x_{i,q,t}^{h_k}=0,\,\forall k \in \mathcal{K}, \, h_k\in\mathcal{H}^k,\, i\in\mathcal{I}^k,\,q\in\mathcal{Q}^k,\, t\in\mathcal{T},\,t+T_{i,q}^{k}+S_{q}\geq T_2
\label{eq:Routing - 5_2}\\
&\sum_{i'\in\mathcal{I}^k}x_{i',j,\left(t-T_{j,i'}^{k}-S_{j}\right)}^{h_k}\geq z_{j,i,t}^{h_k}, \forall k \in \mathcal{K}, \,  h_k\in\mathcal{H}^k,\, i\in\mathcal{I}^k,\, j\in\mathcal{J},\, t\in\mathcal{T},\,t-T_{j,i}^{k}-S_{j}\geq T_1\label{eq:Routing - 10}\\
&\sum_{q\in\mathcal{Q}^k}x_{i,q,\left(t+T_{j,i}^{k}+S_{i}\right)}^{h_k}+r_{i,0,\left(t+T_{j,i}^{k}+S_{i}\right)}^{h_k}\geq z_{j,i,t}^{h_k},\label{eq:Routing - 7}\\
&\forall k \in \mathcal{K}, \, h_k\in\mathcal{H}^k,\, t\in\mathcal{T},\, i\in\mathcal{I}^k,\, j\in\mathcal{J},\,\,t+T_{j,i}^{k}+S_{i}<T_2\nonumber\\
&z_{j,i,t}^{h_k}=0, \forall k \in \mathcal{K}, \,  h_k\in\mathcal{H}^k,\, i\in\mathcal{I}^k,\, j\in\mathcal{J},\, t\in\mathcal{T},\,\,t+T_{j,i}^{k}+S_{i}\geq T_2
\label{eq:Routing - 7_2}\\
&\sum_{j^{'}\in\mathcal{J}-\{j\}}\sum_{\tau=t-S_{j}}^{t+T_{j,i}^{k}+S_{i}}\sum_{i^{'}\in\mathcal{I}^k}z_{j^{'},i^{'},\tau}^{h_k} + \sum_{q^{'}\in\mathcal{Q}^k}\sum_{\tau=t-S_{j}}^{t+T_{j,i}^{k}+S_{i}}\sum_{i^{'}\in\mathcal{I}^k}x_{i^{'},q^{'},\tau}^{h_k}\leq M_{1}\left(1-z_{j,i,t}^{h_k}\right),\label{eq:Routing - 8}\\
&\forall k \in \mathcal{K}, \, h_k\in\mathcal{H}^k,\,i\in\mathcal{I}^k,\, j\in\mathcal{J},\, t\in\mathcal{T},\,t+T_{j,i}^{k}+S_{i}<T_2\nonumber
\end{align}

Constraints (\ref{eq:Routing - 1}) ensure that each ESB $h \in \mathcal{H}$ departs the depot at most once for only one shelter. 
Constraints (\ref{eq:Routing - 2}) enforce that each ESB $h$ must be at most in one of the following states: (1) departing a CS $j$ for a shelter $i$, (2) departing a shelter $i$ for a CS or another shelter $q$, (3) actively serving (discharging) a shelter $i$, or (4) returning to the depot if the ESB $h$ has already departed the depot in the previous time slots. Constraints (\ref{eq:Routing - 3}) ensure that if ESB $h$ departs the depot for shelter $i$ at time slot $t$, it should depart the shelter $i$ at time slot $t+T_{0,i}^k+S_i$. However, if serving shelter $i$ starts after time slot $T_2$, departing the depot is not allowed by constraints (\ref{eq:Routing - 3-1}), as all ESB activities must be completed within the planning horizon. Constraints (\ref{eq:Routing 4 - 1}) ensure that if ESB $h$ departs shelter $i$ for CS $j$ at time slot $t$, the ESB $h$ should depart the CS $j$ at time slot $t+T_{i,j}^k+S_j$. 
Constraints (\ref{eq:Routings 4 - 2}) ensure that if ESB $h$ departs shelter $i$ for a different shelter $i'$, the ESB $h$ should depart shelter $i'$ at time slot $t+T_{i,i'}^k+S_{i'}$. As specified in constraints (\ref{eq:Routing - 5}), if ESB $h$ departs shelter $i$ for another shelter or a CS $q$ at time slot $t$, the ESB $h$ cannot depart from any other shelters or CSs between time slots $t-S_i$ and $t+T_{i,q}^k+S_{q}$. Here, the duration between time slots $t-S_i$ and $t+T_{i,q}^k+S_{q}$, includes the service time of shelter $i$, the travel time from shelter $i$ to another location $q$, and the service time of $q$. If $t+T_{i,q}^k+S_{q}>T_2$, departing shelter $i$ for any node $q$ is not allowed, captured by constraints (\ref{eq:Routing - 5_2}). Constraints (\ref{eq:Routing - 10}) and (\ref{eq:Routing - 7}) represent that if ESB $h$ departs CS $j$ for shelter $i$ at time slot $t$, the ESB $h$ has already departed a shelter $i'$ previously for the CS $j$ at time slot $t-T_{j,i'}^k-S_{j}$, and should depart shelter $i$ at time slot $t+T_{j,i}^k+S_{i}$. However, if departing shelter $i$ happens after the planning horizon, departing CS $j$ for shelter $i$ at time slot $t$ is not allowed by constraints (\ref{eq:Routing - 7_2}). Constraints (\ref{eq:Routing - 8}) ensure that if ESB $h$ departs CS $j$ for shelter $i$ at time slot $t$, between time slots $t-S_j$ and $t+T_{j,i}^k+S_{i}$, ESB $h$ cannot depart any other CSs or shelters. Here, $M_{1}$ needs to be greater than $\max_{k\in\mathcal K} \{C_{max}^{k}\}$.

\paragraph*{Status of ESBs in Serving Shelters:}

\allowdisplaybreaks
\begin{align}
& w_{i,\tau}^{h_k}\geq\sum_{q\in\mathcal{Q}^k}x_{i,q,t}^{h_k}, \forall k \in \mathcal{K}, h_k\in\mathcal{H}^k,\,i\in\mathcal{I}^k,\, t\in\mathcal{T},\,\tau\in\left[t-S_{i},t-1\right]\label{eq:Waiting - 1}\\
& w_{i,\tau}^{h_k}\leq1-x_{i,q,t}^{h_k}, \forall k \in \mathcal{K}, h_k\in\mathcal{H}^k,\, i\in\mathcal{I}^k,\, q\in\mathcal{Q}^k,\, t\in\mathcal{T},\,\tau\in\left[t,t+T_{i,q}^{k}+S_{q}\right]\label{eq:Waiting - 2}\\
& w_{i,\tau}^{h_k}\leq1-z_{j,i,t}^{h_k}, \forall k \in \mathcal{K}, h_k\in\mathcal{H}^k,\, i\in\mathcal{I}^k,\, j\in\mathcal{J},\, t\in\mathcal{T},\,\tau\in\left[t,t+T_{j,i}^{k}-1\right]\label{eq:Waiting - 3}\\
& w_{i,\tau}^{h_k}\geq z_{j,i,t}^{h_k}, \forall k \in \mathcal{K}, h_k\in\mathcal{H}^k,\, i\in\mathcal{I}^k,\, j\in\mathcal{J},\, t\in\mathcal{T},\,\tau\in\left[t+T_{j,i}^{k},t+T_{j,i}^{k}+S_{i}-1\right]\label{eq:Waiting - 4}\\
& w_{i,\tau}^{h_k}\leq1-u_{0,i,t}^{h_k}, \forall k \in \mathcal{K}, h_k\in\mathcal{H}^k,\, i\in\mathcal{I}^k,\, t\in\mathcal{T},\,\tau\in\left[T_1,t+T_{0,i}^{k}-1\right]\label{eq:Waiting - 5}\\
& w_{i,\tau}^{h_k}\geq u_{0,i,t}^{h_k}, \forall k \in \mathcal{K}, h_k\in\mathcal{H}^k,\, i\in\mathcal{I}^k,\, t\in\mathcal{T},\,\tau\in\left[t+T_{0,i}^{k},t+T_{0,i}^{k}+S_{i}-1\right]\label{eq:Waiting - 6}\\
& w_{i^{'},t}^{h_k}\leq1-w_{i,t}^{h_k}, \forall k \in \mathcal{K}, h_k\in\mathcal{H}^k,\, i\in\mathcal{I}^k,\, t\in\mathcal{T},\, i'\in\mathcal{I}^k\setminus\left\{ i\right\} \label{eq:Waiting - 7}\\
&\sum_{\tau=t+1}^{t+S_{i}}\sum_{q\in\mathcal{Q}^k}x_{i,q,t}^{h_k}+\sum_{\tau=t+1}^{t+S_{i}}r_{i,0,\tau}^{h_k}\geq w_{i,t}^{h_k},\forall k \in \mathcal{K}, h_k\in\mathcal{H}^k,\, i\in\mathcal{I}^k,\, t\in\mathcal{T},\,t+S_{i}<T_2\label{eq:Waiting - 8}\\
&\sum_{i^{'}\in\mathcal{I}^k\setminus\left\{i\right\}}\sum_{\tau=t-T_{i^{'},i}^{k}-S_{i}}^{t-T_{i^{'},i}^{k}}x_{i^{'},i,\tau}^{h_k}+\sum_{j\in\mathcal{J}}\sum_{\tau=t-T_{j,i}^{k}-S_{i}}^{t-T_{j,i}^{k}}z_{j,i,\tau}^{h_k}+\sum_{\tau=t-T_{0,i}^{k}-S_{i}}^{t-T_{0,i}^{k}}u_{i,\tau}^{h_k}\geq w_{i,t}^{h_k},\label{eq:Waiting - 8_2}\\
& \forall k \in \mathcal{K}, h_k\in\mathcal{H}^k,\, i\in\mathcal{I}^k,\, t\in\mathcal{T},\,t+S_{i}\geq T_2\nonumber
\end{align}

Constraints (\ref{eq:Waiting - 1}) and (\ref{eq:Waiting - 2}) enforce that if ESB $h$ departs shelter $i$ for another shelter or a CS $q$ at time slot $t$, this ESB should have been actively serving (i.e., discharging energy) shelter $i$ between time slots $t-S_i$ and $t-1$, whereas it cannot serve shelter $i$ from time slots $t$ to $t+T_{i,q}^k+S_q$. Similarly, constraints (\ref{eq:Waiting - 3}) and (\ref{eq:Waiting - 4}) ensure that if ESB $h$ departs CS $j$ for shelter $i$ at time slot $t$, starting from time slot $t$ until time slot $t+T_{j,i}^k-1$, the ESB $h$ cannot be on the status of actively serving shelter $i$, whereas it should actively serve shelter $i$ from time slots $t+T_{j,i}^k$ to $t+T_{j,i}^k+S_i-1$. Constraints (\ref{eq:Waiting - 5}) and (\ref{eq:Waiting - 6}) represent that if ESB $h$ departs the depot for shelter $i$ at time slot $t$, the ESB $h$ should actively serve shelter $i$ from time slots $t+T_{0,i}^k$ to $t+T_{0,i}^k+S_i-1$, whereas it cannot actively serve shelter $i$ until time slot $t+T_{0,i}^k-1$. Constraints (\ref{eq:Waiting - 7}) enforce the rule that if ESB $h$ is actively serving shelter $i$ at time slot $t$, this ESB cannot serve another shelter $i'$ at time slot $t$. Furthermore, constraints (\ref{eq:Waiting - 8}) and (\ref{eq:Waiting - 8_2}) ensure that if ESB $h$ is actively serving shelter $i$ at time slot $t$, this ESB should have been on the status of arriving to shelter $i$ in one of the previous time slots, and should depart shelter $i$ at a time slot between $t+1$ and $t+S_i$.

\paragraph*{Status of ESBs in Returning to the Depot:}

\allowdisplaybreaks
\begin{align}
& \sum_{i\in\mathcal{I}^k}\sum_{\tau=T_1}^{t}r_{i,0,\tau}^{h_k}\leq\sum_{i\in\mathcal{I}^k}\sum_{\tau=T_1}^{t}u_{0,i,\tau}^{h_k}, \forall k \in \mathcal{K}, \, h_k\in\mathcal{H}^k,\, t\in\mathcal{T}\label{eq:Routing - 1-2}\\
& \sum_{i\in\mathcal{I}^k}\sum_{\tau=t+1}^{T_2}r_{i,0,\tau}^{h_k}\geq\sum_{i\in\mathcal{I}^k}u_{0,i,t}^{h_k}, \forall k \in \mathcal{K}, \,  h_k\in\mathcal{H}^k,\, t\in\mathcal{T}\label{eq:Routing - 1-3}\\
& \sum_{\tau=t}^{T_2}\sum_{i\in\mathcal{I}^k}\left(w_{i,\tau}^{h_k}+\sum_{j\in\mathcal{J}}z_{j,i,\tau}^{h_k}+\sum_{q\in\mathcal{Q}^k}x_{i,q,\tau}^{h_k}\right)\leq M_{1}\left(1-\sum_{i\in\mathcal{I}^k}r_{i,0,t}^{h_k}\right), \forall k \in \mathcal{K}, \,  h_k\in\mathcal{H}^k,\, t\in\mathcal{T}\label{eq:Routing - 2-1}\\
& w_{i,\tau}^{h_k}\geq r_{i,0,t}^{h_k}, h_k\in\mathcal{H}^k,\, i\in\mathcal{I}^k,\, t\in\mathcal{T},\,\tau\in\left[t-S_{i},t-1\right]\label{eq:Return - 1}\\
& u_{0,i,T_2}^{h_k},x_{i,q,T_2}^{h_k},z_{j,i,T_2}^{h_k},w_{i,T_2}^{h_k}=0, \forall k \in \mathcal{K}, \, h_k\in\mathcal{H}^k,\, i\in\mathcal{I}^k,\, q \in \mathcal{Q}^k ,\, j \in \mathcal{J}\label{eq:Return - 1-1}
\end{align}

Constraints (\ref{eq:Routing - 1-2}) and (\ref{eq:Routing - 1-3}) ensure that an ESB returns to the depot, if and only if, it departs the depot as well. Moreover, constraints \eqref{eq:Routing - 1-3} enforce that returning to the depot occurs after departing the depot for a shelter. As specified by constraints (\ref{eq:Routing - 2-1}), if ESB $h$ returns to the depot at time slot $t$, in the other time slots after this time, this ESB cannot be in any other possible states, including departing a shelter, actively serving a shelter, or departing a CS. Constraints (\ref{eq:Return - 1}) ensure that if ESB $h$ returns to the depot from shelter $i$ at time slot $t$, the ESB $h$ should have been actively serving (i.e., discharging energy) shelter $i$ between time slots $t-S_i$ and $t-1$. As we assume that all ESBs should finally return to the depot, no activity except for returning to the depot is allowed at the last time slot, which is represented by constraints (\ref{eq:Return - 1-1}).

\paragraph*{Energy Scheduling Constraints:}

\allowdisplaybreaks
\begin{align}
& C_{max}^{k}-R_{0,i}^{k}-M_{1}*\left(1-u_{0,i,t}^{h_k}\right)\leq soc_{\left(t+T_{0,i}^{k}\right)}^{h_k}\leq C_{max}^{k}-R_{0,i}^{k}+M_{1}*\left(1-u_{0,i,t}^{h_k}\right), \label{eq:Energy - 4}\\
& \forall k \in \mathcal{K}, \,  h_k\in\mathcal{H}^k,\, i\in\mathcal{I}^k,\, t\in\mathcal{T}\nonumber\\
& soc_{t}^{h_k}-R_{i,q}^{k}-M_{1}*\left(1-x_{i,q,t}^{h_k}\right)\leq soc_{\left(t+T_{i,q}^{k}\right)}^{h_k}\leq soc_{t}^{h_k}-R_{i,q}^{k}+M_{1}*\left(1-x_{i,q,t}^{h_k}\right)\label{eq:Energy - 5}\\
& ,\forall k \in \mathcal{K}, \,  h_k\in\mathcal{H}^k,\, i\in\mathcal{I}^k,\, q\in\mathcal{Q}^k,\, t\in\mathcal{T}\nonumber\\
& C_{max}^{k}-R_{j,i}^{k}-M_{1}*\left(1-z_{j,i,t}^{h_k}\right)\leq soc_{\left(t+T_{j,i}^{k}\right)}^{h_k}\leq C_{max}^{k}-R_{j,i}^{k}+M_{1}*\left(1-z_{j,i,t}^{h_k}\right),\label{eq:Energy - 6}\\
& \forall k \in \mathcal{K}, \,  h_k\in\mathcal{H}^k,\, i\in\mathcal{I}^k,\, j\in\mathcal{J},\, t\in\mathcal{T}\nonumber \\
& -M_{2}*\left(1-\sum_{q\in\mathcal{Q}^k}x_{i,q,t}^{h_k}\right)\leq soc_{t}^{h_k}-soc_{\left(t-S_{i}\right)}^{h_k}+\sum_{\tau=t-S_{i}}^{t-1}\left(g_{i,\tau}^{h_k}\right)\leq M_{2}*\left(1-\sum_{q\in\mathcal{Q}^k}x_{i,q,t}^{h_k}\right),\label{eq:Energy - 7}\\
& \forall k \in \mathcal{K}, \,  h_k\in\mathcal{H}^k,\, i\in\mathcal{I}^k,\, t\in\mathcal{T}\nonumber \\
& -M_{2}*\left(1-r_{i,0,t}^{h_k}\right)\leq soc_{t}^{h_k}-soc_{\left(t-S_{i}\right)}^{h_k}+\sum_{\tau=t-S_{i}}^{t-1}\left(g_{i,\tau}^{h_k}\right)\leq M_{2}*\left(1-r_{i,0,t}^{h_k}\right),\label{eq:Energy - 7-1}\\
& \forall k \in \mathcal{K}, \, h_k\in\mathcal{H}^k,\, i\in\mathcal{I}^k,\, t\in\mathcal{T}\nonumber \\
& soc_{t}^{h_k}\geq R_{i,0}^{k} r_{i,0,t}^{h_k}+C_{min}^{k},\forall k \in \mathcal{K},\, h_k\in\mathcal{H}^k,\, i\in\mathcal{I}^k,\, t\in\mathcal{T}\label{eq:Energy - 7-2}\\
& C_{min}^{k}\leq soc_{t}^{h_k}\leq C_{max}^{k}, \forall k \in \mathcal{K}, \, h_k\in\mathcal{H}^k,\, t\in\mathcal{T}\label{eq:Energy - 8}
\end{align}

Constraints (\ref{eq:Energy - 4}) compute the SOC of ESB $h$ at the arrival time to shelter $i$ after departing the depot with a fully charged battery. Similarly, constraints (\ref{eq:Energy - 5}) compute the SOC of ESBs upon arrival at a shelter or a CS $q$ from another shelter, considering the energy consumption on that route (i.e., $R_{i,q}^{k}$). Constraints (\ref{eq:Energy - 6}) calculate the SOC of ESBs upon arrival at a shelter from a CS with a fully charged battery. Moreover, constraints (\ref{eq:Energy - 7}) and (\ref{eq:Energy - 7-1}) compute the SOC of ESB $h$ after serving shelter $i$. If ESB $h$ departs shelter $i$ for another shelter or CS, or returns to the depot at time slot $t$, the ESB $h$ has arrived at this shelter $i$ at a previous time slot $t-S_{i}$. Therefore, the difference between the SOC of ESB $h$ at time slots $t$ and $t-S_{i}$ equals the amount of discharged energy from ESB $h$ to shelter $i$ during the service time of shelter $i$ (i.e., $S_i$). Constraints (\ref{eq:Energy - 7-2}) ensure that each ESB $h$ has sufficient SOC before returning to the depot. Finally, constraints (\ref{eq:Energy - 8}) ensure that the SOC of each ESB remains within the acceptable range at all times. Here, $M_{2}$ must be greater than $\sum_{i\in\mathcal{I}}\sum_{t\in\mathcal{T}}P_{i,t}$.

\paragraph*{Relationship Between ESB Routing and Energy Scheduling:}

\begin{equation}
G_{min}^{k}w_{i,t}^{h_k}\leq g_{i,t}^{h_k}\leq M_{1}*w_{i,t}^{h_k}, \forall k \in \mathcal{K}, \,  h_k\in\mathcal{H}^k,\, i\in\mathcal{I}^k,\, t\in\mathcal{T}\label{eq:Energy - 9}
\end{equation}

Constraints (\ref{eq:Energy - 9}) ensure that energy discharge by an ESB to a shelter is possible, if and only if, the ESB actively serves the shelter (i.e., $w_{i,t}=1$). We enforce $g_{i,t}^{h_k} \geq G_{min}^{k}$ when an ESB actively serves a shelter.

\paragraph*{Sign Restrictions of Decision Variables:}
\allowdisplaybreaks
\begin{align}
    &x_{i,q,t}^{h_k}, z_{j,i,t}^{h_k}, w_{i,t}^{h_k}, u_{0,i,t}^{h_k}, r_{i,0,t}^{h_k} \in \{0, 1\};\, l_{i}, g_{i,t}^{h_k}, soc_{t}^{h_k} \in \mathbb{R}^{+}, \, g_{i,t}^{h_k} \leq C_{max}^k - C_{min}^k\label{eq:Continuous_Variables}
\end{align}

Constraints (\ref{eq:Continuous_Variables}) define the restrictions on the binary and continuous decision variables.

\paragraph*{Valid Inequalities:}
As we have multiple available ESBs of the same type that bring symmetry into the problem, we introduce the symmetry-breaking constraints as valid inequalities to improve the computational efficiency of the \textit{ESB-BaM-MILP} (i.e., Eqs. \eqref{eq:Obj Function} - \eqref{eq:Continuous_Variables}).
Constraints (\ref{eq:Symmetry Breaking}) ensure that for each ESB type $k$, if an ESB $h_k$ does not participate in the restoration program, no ESB of the same type with a higher index than $h_k$ can participate in restoration. Constraints (\ref{eq:Symmetry Breaking-1}) are tighter versions of constraints (\ref{eq:Symmetry Breaking}), defined for each time slot $t$. Therefore, we only use constraints (\ref{eq:Symmetry Breaking-1}) in our numerical experiments.     
\allowdisplaybreaks
\begin{align}
& \sum_{h_k^{'}=h_k+1}^{H_{max}^k}\sum_{i\in\mathcal{I}^k}\sum_{t\in\mathcal{T}}u_{0,i,t}^{h_k^{'}}\leq M_{1}\sum_{i\in\mathcal{I}^k}\sum_{t\in\mathcal{T}}u_{0,i,t}^{h_k},\qquad\forall h_k\in\mathcal{H}^k, \forall k \in \mathcal{K}\label{eq:Symmetry Breaking}\\
& \sum_{h_k^{'}=h_k+1}^{H_{max}^k}\sum_{i\in\mathcal{I}^k}\sum_{\tau=0}^{t}u_{0,i,\tau}^{h_k^{'}}\leq M_{1}\sum_{i\in\mathcal{I}^k}\sum_{\tau=0}^{t}u_{0,i,\tau}^{h_k},\qquad\forall h_k\in\mathcal{H}^k,\,\forall t\in\mathcal{T}, \forall k \in \mathcal{K}\label{eq:Symmetry Breaking-1}
\end{align}

\begin{remark}
    The symmetry-breaking constraints enforce that if an ESB $h_k$ of type $k$ is not dispatched, then a higher-indexed ESB such that $h_k' > h_k$ of the same type $k$ must also not be dispatched, eliminating redundant (i.e., symmetric) solutions. This is because, if an ESB $h_k'$ of type $k$ is used in an optimal solution, then any ESB of the same type with a lower index $h_k < h_k'$ can also replace it without changing the objective value.
\end{remark}

\begin{prop}\label{prop:symmetry_breaking}
    The \textit{ESB-BaM-MILP} model augmented with symmetry-breaking constraints (\ref{eq:Symmetry Breaking}) - (\ref{eq:Symmetry Breaking-1}) remains as a valid model, and its feasible set is a proper subset of that of \textit{ESB-BaM-MILP}.
\end{prop}

Indeed, we mention that \textit{ESBILR}, modeled as \textit{ESB-BaM-MILP}, is a theoretically challenging problem, which generalizes the classical Traveling Salesman Problem (TSP).
\begin{prop}\label{prop:NP-hardness}
    The \textit{ESBILR} is NP-hard.
\end{prop}

We define the model presented in Eqs. (\ref{eq:Obj Function}) - (\ref{eq:Continuous_Variables}) and (\ref{eq:Symmetry Breaking-1}) as the accelerated MILP model, and is referred to as \textit{ESB-AccM-MILP}.
The proofs of Propositions \ref{prop:symmetry_breaking} and \ref{prop:NP-hardness} are provided in Appendix \ref{Appendix_Proof}.

\section{Branch and Price Algorithm}\label{sec:BnP_Algorithm}

In this section, we present an improved B\&P algorithm that combines column generation (CG) and branch-and-bound (B\&B) techniques to efficiently solve the \textit{ESBILR}. We first reformulate the \textit{ESB-BaM-MILP} (i.e., Eqs. (\ref{eq:Obj Function}) - (\ref{eq:Continuous_Variables})) to a new route-based set covering formulation presented in model (\ref{eq:Set_Covering_Model}) of Section \ref{subsec:Set_Covering_Formulation}. We then use Dantzig-Wolfe decomposition to decompose the set covering model (\ref{eq:Set_Covering_Model}) into a master problem and pricing subproblems. In the context of our problem, the subproblem generates new routes (i.e.,  columns), where a \textit{route} represents transportation and energy scheduling decisions for ESBs. The master problem is responsible for assigning ESBs to routes to satisfy the energy demand of shelters with a minimum cost. As enumerating all feasible routes is computationally impractical, we employ the CG scheme at each node of the B\&B tree to efficiently solve that node, which involves iteratively solving a LP-relaxed master problem (\textit{RLMP}) with a restricted set of initial feasible routes, $\mathcal{P}_{k}$, indexed by $p_{k}$, and pricing subproblems. This section also presents a heuristic B\&P variant that leverages a labeling-integrated dynamic programming algorithm to solve the pricing subproblems faster. 

\subsection{Set Covering Formulation}\label{subsec:Set_Covering_Formulation}

Model (\ref{eq:Set_Covering_Model}) presents the set covering formulation, where $\mathcal{P}_{k}^{'}$ is the set of all feasible routes for ESB type $k$. As mentioned above, each route contains transportation and energy scheduling decisions for ESBs. Moreover, in our problem, multiple ESBs can follow identical routes. Therefore, we define $\lambda_{p_{k}}^{k}$ as an integer variable representing the number of ESBs of type $k$ using route $p_k$. 
\begin{subequations} \label{eq:Set_Covering_Model}
\begin{align}
\min \, & \sum_{k\in\mathcal{K}}\sum_{p_{k}\in\mathcal{P}_{k}^{'}}C_{p_{k}}^{k}\lambda_{p_{k}}^{k} + \sum_{i\in\mathcal{I}}F_{i}l_{i}\label{eq:SCa}\\
\text{s.t.} \, & \sum_{k\in\mathcal{K}^i}\sum_{t\in\mathcal{T}}\sum_{p_{k}\in\mathcal{P}_{k}^{'}}G_{i,t,p_{k}}\lambda_{p_{k}}^{k} + l_{i} \geq \sum_{t\in\mathcal{T}}P_{i,t}, \quad \forall i \in \mathcal{I} \label{eq:SCb}\\
& l_i \geq 0 ,\, \lambda_{p_k}^{k} \in \mathbb{Z} \nonumber
\end{align}
\end{subequations}

The objective function (\ref{eq:MPa}) is equivalent to the objective function of the \textit{ESB-BaM-MILP} (i.e., Eq. (\ref{eq:Obj Function})), which aims to minimize the total investment and operational costs, as well as the penalty cost of aggregated unmet demand. $C_{p_{k}}^{k}$ is the total investment and energy consumption cost of route $p_k$. Constraints (\ref{eq:MPb}) ensure that the aggregated demand of each shelter $i$ over the planning horizon (i.e., $\sum_{t\in\mathcal{T}}P_{i,t}$) is satisfied either fully by the ESBs on their respective routes (where $G_{i,t,p_{k}}$ is the amount of discharged energy of the ESB to shelter $i$ at time slot $t$ in route $p_k$), or by incurring penalty on unmet demand ($l_i$) when the ESBs cannot satisfy the entire demand.

\subsection{Master Problem}\label{subsubsec:Master_Problem}

The mathematical formulation for the \textit{RLMP} is as the model (\ref{eq:master_problem}).
\begin{subequations} \label{eq:master_problem} 
\begin{align}
\min \, & \sum_{k\in\mathcal{K}}\sum_{p_{k}\in\mathcal{P}_{k}}C_{p_{k}}^{k}\lambda_{p_{k}}^{k} + \sum_{i\in\mathcal{I}}F_{i}l_{i}\label{eq:MPa}\\
\text{s.t.} \, & \sum_{k\in\mathcal{K}^i}\sum_{t\in\mathcal{T}}\sum_{p_{k}\in\mathcal{P}_{k}}G_{i,t,p_{k}}\lambda_{p_{k}}^{k} + l_{i} \geq \sum_{t\in\mathcal{T}}P_{i,t}, \quad \forall i \in \mathcal{I} 
 & (\pi_i) \label{eq:MPb}\\
& n_{EB}^{k} = \sum_{p_{k}\in\mathcal{P}_{k}}\lambda_{p_{k}}^{k}, \quad \forall k \in \mathcal{K} & (\psi^k) \label{eq:MPd}\\
& f_i = \sum_{t \in \mathcal{T}} \sum_{k \in \mathcal{K}^i} \sum_{p_{k} \in \mathcal{P}_k} W_{i,t,p_k} \lambda_{p_k}^{k}, \quad \forall i \in \mathcal{I} & (\mu_i) \label{eq:MPe}\\
& b_{i}^{k} = \sum_{t \in \mathcal{T}} \sum_{p_{k} \in \mathcal{P}_k} W_{i,t,p_k} \lambda_{p_k}^{k}, \quad \forall i \in \mathcal{I}, \forall k \in \mathcal{K}^i & (\mu_i^k) \label{eq:MPf}\\
& d_{i,t} = \sum_{k \in \mathcal{K}^i} \sum_{p_{k} \in \mathcal{P}_k} W_{i,t,p_k} \lambda_{p_k}^{k}, \quad \forall i \in \mathcal{I}, \forall t \in \mathcal{T} & (\rho_{i,t}) \label{eq:MPg}\\
& m_{i,t}^{k} = \sum_{p_{k} \in \mathcal{P}_k} W_{i,t,p_k} \lambda_{p_k}^{k}, \quad \forall i \in \mathcal{I}, \forall t \in \mathcal{T}, \forall k \in \mathcal{K}^i & (\rho_{i,t}^k) \label{eq:MPh}\\
& v_{i,j} = \sum_{t \in \mathcal{T}} \sum_{k \in \mathcal{K}^i} \sum_{p_k \in \mathcal{P}_k} \left(X_{i,j,t,p_k} + Z_{j,i,t,p_k}\right) \lambda_{p_k}^{k}, \quad \forall i \in \mathcal{I}, \forall j \in \mathcal{J} & (\eta_{i,j}) \label{eq:MPi}\\
& w_{i,j}^{k} = \sum_{t \in \mathcal{T}} \sum_{p_k \in \mathcal{P}_k} \left(X_{i,j,t,p_k} + Z_{j,i,t,p_k}\right) \lambda_{p_k}^{k}, \quad \forall i \in \mathcal{I}, \forall j \in \mathcal{J}, \forall k \in \mathcal{K}^i & (\eta_{i,j}^k) \label{eq:MPj}\\
& a_{i} = \sum_{t \in \mathcal{T}} \sum_{k \in \mathcal{K}^i} \sum_{p_k \in \mathcal{P}_k} \left(U_{0,i,t,p_k} + R_{i,0,t,p_k}\right) \lambda_{p_k}^{k}, \quad \forall i \in \mathcal{I}, \forall j \in \mathcal{J} & (\theta_i) \label{eq:MPk}\\
& e_{i}^{k} = \sum_{t \in \mathcal{T}} \sum_{p_k \in \mathcal{P}_k} \left(U_{0,i,t,p_k} + R_{i,0,t,p_k}\right) \lambda_{p_k}^{k}, \quad \forall i \in \mathcal{I}, \forall j \in \mathcal{J}, \forall k \in \mathcal{K}^i & (\theta_i^k) \label{eq:MPl}\\
& \lambda_{p_k}^{k}, n_{EB}^{k}, f_i, d_{i,t}, v_{i,j}, w_{i,j}^{k}, a_i, e_{i}^{k}, b_{i}^{k}, m_{i,t}^{k}, l_i \geq 0  \nonumber
\end{align}
\end{subequations}

The objective function of the \textit{RLMP}, defined in Eq. (\ref{eq:MPa}), and constraints (\ref{eq:MPb}) are equivalent to the model (\ref{eq:Set_Covering_Model}), with the differences that the \textit{RLMP} works with a restricted set of feasible routes (i.e., $\mathcal{P}_{k}$), and the LP-relaxation of the integer $\lambda_{p_{k}}^{k}$ variables. In Eqs. \eqref{eq:MPb} - \eqref{eq:MPl}, the Greek letters inside parentheses beside each constraint (e.g., $\pi_i$) represent the dual variables associated with the corresponding constraints.

At each node in the B\&B tree, the \textit{RLMP} (\ref{eq:master_problem}) may not necessarily converge to an integral solution with respect to $\lambda_{p_k}^k$. Therefore, we apply different branching rules to progressively enforce integrality as the algorithm proceeds through the B\&B tree. 
We introduce branching rules as equality constraints (\ref{eq:MPd}) - (\ref{eq:MPl}) that are hierarchically applied in the following order: (1) branching on $n_{EB}^k$ (the number of utilized ESBs of type $k$); (2) branching on $f_i$, $d_{i,t}$, $b_i^k$, and $m_{i,t}^k$ (the number of total and timely visits to shelter $i$ by all ESBs and ESBs of type $k$, respectively); (3) branching on $v_{i,j}$, $w_{i,j}^k$, $a_{i}$, and $e_{i}^k$ (the number of times the arc $i-j$ or $i-0$ is traversed over the planning horizon, respectively); and (4) branching explicitly on the fractional $\lambda_{p_k}^k$ variables, if none of the previous branching rules are applicable but the solution is not integer.
Here, $W_{i,t,p_{k}}$ is a binary parameter that represents whether on route $p_k$ the ESB serves shelter $i$ at time slot $t$. In Eqs. (\ref{eq:MPi}) - (\ref{eq:MPj}), $X_{i,j,t,p_{k}}$ is a binary parameter that represents whether on route $p_k$ the ESB departs shelter $i$ for CS $j$ at time slot $t$. Similarly, binary parameter $Z_{j,i,t,p_{k}}$ represents whether the ESB departs CS $j$ for shelter $i$ in the same route and time slot. Furthermore, In Eqs. (\ref{eq:MPk}) - (\ref{eq:MPl}), the parameter $U_{0,i,t,p_{k}}$ is a binary parameter that represents whether on route $p_k$ the ESB departs the depot for shelter $i$ at time slot $t$, whereas $R_{i,0,t,p_{k}}$ represents whether on route $p_k$ the ESB returns to the depot from shelter $i$ at time slot $t$.

\subsection{Exact Pricing Subproblems}\label{subsubsec:Pricing_Subproblems}

As we have different ESB types with different features---battery capacities, energy consumption rates, and costs---the route for each ESB type can be different. Therefore, we need to define different subproblems, one for each ESB type. A column (route) generated by the subproblem must be feasible with respect to the routing and energy scheduling constraints of the \textit{ESB-BaM-MILP} (i.e., Eqs. (\ref{eq:Obj Function}) - (\ref{eq:Continuous_Variables})).
The mathematical formulation of the exact subproblem for generating a route for ESBs of type $k$ is presented in model (\ref{eq:Subproblem}).
\begin{align}
\label{eq:Subproblem}
& min \qquad \Bigg\{ C_{fx}\sum_{i\in\mathcal{I}^k}\sum_{t\in\mathcal{T}}\left(\sum_{q\in\mathcal{Q}}\left(x_{i,q,t}R_{i,q}^{k}\right) + \sum_{j\in\mathcal{J}}\left(z_{j,i,t}R_{j,i}^{k}\right) + \left(u_{0,i,t}R_{0,i}^{k}\right) + \left(r_{i,0,t}R_{i,0}^{k}\right)\right)\\\nonumber
 & \qquad \qquad + \sum_{t\in\mathcal{T}}\sum_{i\in\mathcal{I}^k}C_{inv}^{k}u_{0,i,t}+\sum_{i\in\mathcal{I}^k}\sum_{t\in\mathcal{T}}w_{i,t}\left\{ \mu_{i}+\mu_{i}^{k}+\rho_{i,t}+\rho_{i,t}^{k}\right\} +\psi^{k}\\ \nonumber
 & \qquad \qquad +\sum_{i\in\mathcal{I}^k}\sum_{j\in\mathcal{J}}\left(x_{i,j,t}+z_{j,i,t}\right)\left\{ \eta_{i,j}+\eta_{i,j}^{k}\right\} +\sum_{i\in\mathcal{I}^k}\left(u_{0,i,t}+r_{i,0,t}\right)\left\{ \theta_{i}+\theta_{i}^{k}\right\} -\sum_{i\in\mathcal{I}}\sum_{t\in\mathcal{T}}\left(\left(\pi_{i}\right)g_{i,t}\right) \Bigg\}\\\nonumber
 & s.t. \qquad (\ref{eq:Routing - 1}) - (\ref{eq:Continuous_Variables})
\end{align}

The objective function of the model (\ref{eq:Subproblem}) aims to determine a new route (i.e., routing and energy scheduling decisions) with the most negative reduced cost. 
This subproblem is constrained by the routing and scheduling constraints presented in Eqs. (\ref{eq:Routing - 1}) - (\ref{eq:Continuous_Variables}) (without ESB index, $h_k$) to ensure feasibility of the newly generated route with respect to the routing and energy scheduling constraints.

\subsection{Fast Labeling-Integrated Dynamic Programming for Subproblem}\label{subsubsubsec:LA_Subproblem}

Various labeling algorithms (LAs) have been designed and implemented in the literature (e.g.,~\cite{sadykov2021bucket,moreno2024crew,yamin2024reliable}) as effective approaches to solve pricing subproblems, particularly in the context of vehicle routing problems (VRPs). The flexibility of LA in handling complex constraints, its ability to efficiently prune suboptimal solutions through dominance rules, and its adaptability to various problem-specific extensions, make LA a highly advantageous method for route generation in subproblems. Therefore, we develop a customized labeling algorithm to compute the subproblem in the model \eqref{eq:Subproblem}. 
Unlike traditional VRP, our subproblem allows ESBs of type $k$ to make multiple back-and-forth trips between shelter nodes (for discharging) and CS nodes (for recharging). The key challenge lies in determining not only the routing sequence but also when and how to utilize the CSs to minimize the reduced cost, considering both shelter restoration (i.e., discharging energy to shelters) and battery recharging (i.e., recharging the battery of ESBs in CSs). To address this challenge, we propose an integrated solution framework comprising a Dynamic Programming (DP) layer and an LA layer, as shown in Figure \ref{fig:enter-label}. 
In the lower level, we consider the routing of an ESB over shelters between a pair of CSs, which can be achieved by designing a typical LA. Then, we introduce DP on top of this LA to determine a complete solution to the model \eqref{eq:Subproblem}. Note that such an integrated design matches the network structure and avoids complicated operations to handle functionally different nodes simultaneously. Actually, it often generates an optimal or almost-optimal solution in an efficient fashion.

\subsubsection{An LA for Routing between Charging Stations}

For a given pair of CSs (including the depot), denoted as $l$ and $l'$, we design a regular LA to derive a route linking them. Such a route has the least cost, noting that each arc/node's cost has been modulated by the appropriate shadow prices (i.e., dual values of the \textit{RLMP} \eqref{eq:master_problem}). Indeed, we can treat it as a prize-collecting route (between CSs), as an ESB visits a node, only if, it helps achieve a smaller reduced cost. Suppose that we have an initial route between $l$ and $l'$. The LA will systematically extend this route over the transportation network while adhering to energy and time constraints (i.e., Eqs. \eqref{eq:Routing - 1} - \eqref{eq:Continuous_Variables}).  

Specifically, let label $L^{l',l,t',t}$ be a multi-dimensional representation of a (partial) route, departing node $l'$ at time slot $t'$ and arriving at node $l$ at time slot $t$, which captures key state variables necessary for evaluating or extending this route. In the following derivation, as $(l',l,t',t)$ is given in the context, we simply suppress this superscript to simplify our exposition. Additionally, let $o(L)$ be the last node in the current partial route, $f(L)$ the arrival time at the last node, $soc(L)$ the current energy status of the associated ESB, $G(L)$ the energy discharge pattern of the (partial) route, and $rc(L)$ its reduced cost. When extending this route (label $L$) along an arc to node $y$, a new label $L'$ for the resulting extended partial route is obtained by applying the update rules described in Eqs. \eqref{L1} - \eqref{L5}:
\vspace{-10pt}
\begin{align}
    & o(L') \leftarrow y \label{L1}\\
    & f(L') \leftarrow f(L) + T_{(o(L),y)}^k + S_{y} \label{L2}\\
    & soc(L') \leftarrow soc(L) - g(y) - R_{((o(L),y)}^k \label{L3}\\
    & G(L') \leftarrow \left(G(L),g(y) \right)\label{L4}\\
    & rc(L') \leftarrow rc(L) + irc{({o(L),y})}\label{L5}
\end{align}

where $g(y)$ is the energy discharge pattern for node $o(L')$ (i.e., $y$), and $irc{({o(L),y})}$ denotes the changing value of the reduced cost of this route, extending from $o(L)$ to $y$.

For the new label $L'$, the feasibility of the associated partial route is verified by checking energy and time constraints at node $y$, noting that the energy status cannot be below the minimum battery capacity, and the arrival time cannot exceed time slot $t$. Infeasible $L'$ will be discarded. The complete pseudocode of the LA is presented in Algorithm \ref{alg-LA}. 

\begin{algorithm}[htbp]
\caption{Labeling Algorithm \label{alg-LA}}
\begin{algorithmic}[1]
	\State {\bf Input:} Source and destination CS/depot $l$ and $l'$ with departure and arrival time slots $t$ and $t'$.
        \State {\bf Initialize:} Set $L_0$ with $o(L_0) = l$, $f(L_0) = t$, $soc(L_0) = C_{max}^k$ and $G(L_0)=()$. Also, set list $\mathbb A = \{L_0\}$ and list $\mathbb B = \emptyset$.
	\While {$\mathbb A\neq \emptyset$} 
	\State Select a label $L$ from the list $\mathbb A$;
            \For{Each outgoing arc from node \textbf{$o(L)$}}
            \State Generate $L'$ through Eqs. \eqref{L1} -\eqref{L5}
            \If{$f(L') > t'$ or $soc(L')<C_{min}^k$ }
            \State Remove this label from $\mathbb A$
            \Else
            \If{ $o(L') == l'$}
            \State $\mathbb B = \mathbb B \cup \{L'\}$ 
            \Else 
            \State $\mathbb A = \mathbb A\backslash\{ L\} \cup\{L'\}$
            \EndIf
            \EndIf
            \EndFor
		\EndWhile
		\State {\bf Output:} The route associated with the label of the smallest reduced cost in $\mathbb B$.
\end{algorithmic}\label{1}%
\end{algorithm}

One critical challenge of LA is that the number of labels grows exponentially. As an effective strategy, dominance rules can reduce feasible labels by eliminating labels that are worse than others. Specific to our labeling algorithm, we introduce two dominance rules specified in Proposition \ref{prop_dominance}. 

\begin{prop}
\label{prop_dominance}
Let $\mathbb A$ denote the set of feasible labels obtained so far, and consider two labels $L$ and $L'$ in $\mathbb A$.
$(i)$ If $o(L)=o(L')$, $rc(L')<rc(L)$,  $soc(L')\geq soc(L),$ and $f(L')<f(L)$, we can discard $\{L\}$ from $\mathbb A$ without losing any optimal route. $(ii)$ If $L$ and $L'$ visit the same set of nodes, $rc(L')<rc(L),$ and $f(L')<f(L)$, we can discard $\{L\}$ from $\mathbb A$ without losing any optimal route.
\end{prop}

The proof of Proposition \ref{prop_dominance} is provided in Appendix \ref{Appendix_Proof}.

\subsubsection{DP for Routing over Charging Stations}

The DP layer focuses on high-level decision-making on CS visits; i.e., it determines the optimal sequence of CS visits, while minimizing the reduced cost, considering transportation energy consumption costs, charging (discharging) decisions, and dual variables of the \textit{RLMP} \eqref{eq:master_problem}. The DP process is formulated as a recursive optimization problem with well-defined states, transitions, and boundary conditions. Each state, denoted as $(l, t)$, encapsulates key information about decisions to have the ESB arriving at CS $l$ at time slot $t$. 
Let $RC(l, t)$ be the reduced cost of reaching $(l,t)$, which can be calculated recursively by evaluating all possible transitions from previous states, as described by Eq. \eqref{eq_DP}:

\begin{equation}
\label{eq_DP}
RC(l, t) = \min_{l{'} \in\mathcal{L}, 0\leq
t'\leq t} \{ RC(l^{'},t') +  \min_{L} \{rc(L^{l^{'},l,t',t})\} \}, \ \forall l \in \mathcal{L}, \forall t \in \mathcal{T}
\end{equation}

where $\mathcal{L} = \mathcal{J} \cup \{0\}$, and $rc(L^{l',l,t',t})$ is the incremental reduced cost associated with a transition from CS $l^{'}$ to $l$, over which the smallest one will be selected to update $RC(l,t)$. Note that $rc(L^{l',l,t',t})$ is derived in Eqs. \eqref{L1} - \eqref{L5}. To ensure computational correctness and efficiency, we define the boundary conditions in Eqs. \eqref{b1} - \eqref{b2}. 
\begin{align}
    & RC(0, t) = 0,\  \forall t \in\mathcal{T} \label{b1}\\
    & RC(j, t) = +\infty, \  \forall j \in \mathcal{J},\  \text{if} \  t < t_j^{min} \  \forall t \in\mathcal{T} \label{b2} 
\end{align}

When an ESB is located at the depot at time slot $t$, the reduced cost is initialized as Eq. \eqref{b1}, indicating the zero cost of starting the route at the depot. Additionally, each CS node has a minimum feasible arrival time, denoted as $t_j^{min}$. If $t < t_j^{min}$, the reduced cost is set to infinity as indicated in Eq. \eqref{b2}. 

\begin{figure}[htbp]
    \centering
    \includegraphics[width=0.95\linewidth]{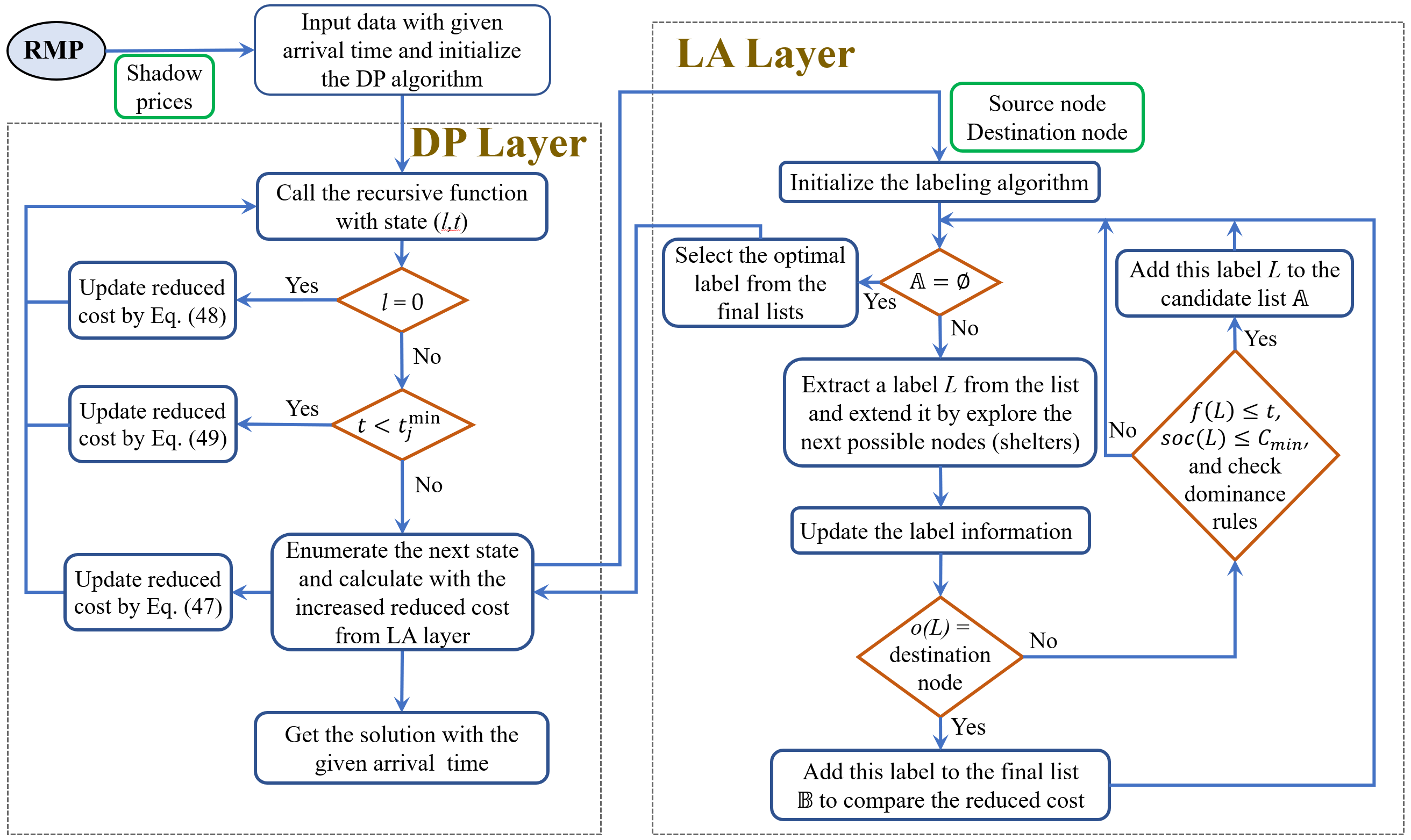}
    \caption{The flowchart of the LA-integrated DP.}
    \label{fig:enter-label}
\end{figure}

As the entire DP procedure is fully defined by Eqs. \eqref{eq_DP} and (\ref{b1}) - (\ref{b2}), we do not provide the detailed steps to simplify our exposition. The flowchart of the LA-integrated DP is presented in Figure \ref{fig:enter-label}.

\subsection{Enhancement Strategies}\label{subsec:Enhancement_Strategies}

In this section, we propose two enhancement strategies to improve the computational efficiency of the B\&P algorithm while maintaining the solution quality.

\paragraph*{Tree Exploration:}
To guide the B\&B tree exploration and improve the lower bound faster, we adopt a \textit{best-first search strategy} that, when selecting the next node to explore, prioritizes the node whose parent node has the smallest objective value among all parent leaf nodes.

\paragraph*{Good Quality Upper Bound:}

After exploring a few number of nodes in the B\&B tree, once a sufficient number of high-quality routes has been generated by CG, we solve the MILP version of the \textit{RLMP} (i.e., the \textit{MILP-MP}) as shown in the model (\ref{eq:MILP_master_problem}). In model (\ref{eq:MILP_master_problem}), by enforcing integrality on $\lambda_{p_k}^k$, the fractional solutions of \textit{RLMP} are converted into integer solutions. The rationale is that, by solving \textit{MILP-MP} after a certain number of nodes in the B\&B tree, the CG has already identified several high-quality feasible routes (i.e., $\mathcal{P}_{k}$), likely to be close to the optimal routes. Therefore, this step generates a feasible integer solution that provides a high-quality upper bound for the original problem, pruning some nodes in the tree, and improving the convergence. The decision on when to solve the \textit{MILP-MP} (\ref{eq:MILP_master_problem}) in the B\&B tree is determined based on the parameter tuning that is discussed in Section \ref{subsec:Experimental_Setup}.
\vspace{-10pt}
\begin{subequations} \label{eq:MILP_master_problem}
\begin{align}
\min \, & \sum_{k\in\mathcal{K}}\sum_{p_{k}\in\mathcal{P}_{k}}C_{p_{k}}^{k}\lambda_{p_{k}}^{k} + \sum_{i\in\mathcal{I}}F_{i}l_{i}\label{eq:MILP-MPa}\\
\text{s.t.} \, & \sum_{k\in\mathcal{K}^i}\sum_{t\in\mathcal{T}}\sum_{p_{k}\in\mathcal{P}_{k}}G_{i,t,p_{k}}\lambda_{p_{k}}^{k} + l_{i} \geq \sum_{t\in\mathcal{T}}P_{i,t}, \quad \forall i \in \mathcal{I} \label{eq:MILP-MPb}\\
& l_i \geq 0 , \, \lambda_{p_k}^{k} \in \mathbb{Z}  \nonumber
\end{align}
\end{subequations}

\section{Case Study}\label{sec:Case_Study}

\begin{figure}[h!]
\begin{centering}
\includegraphics[scale=0.28]{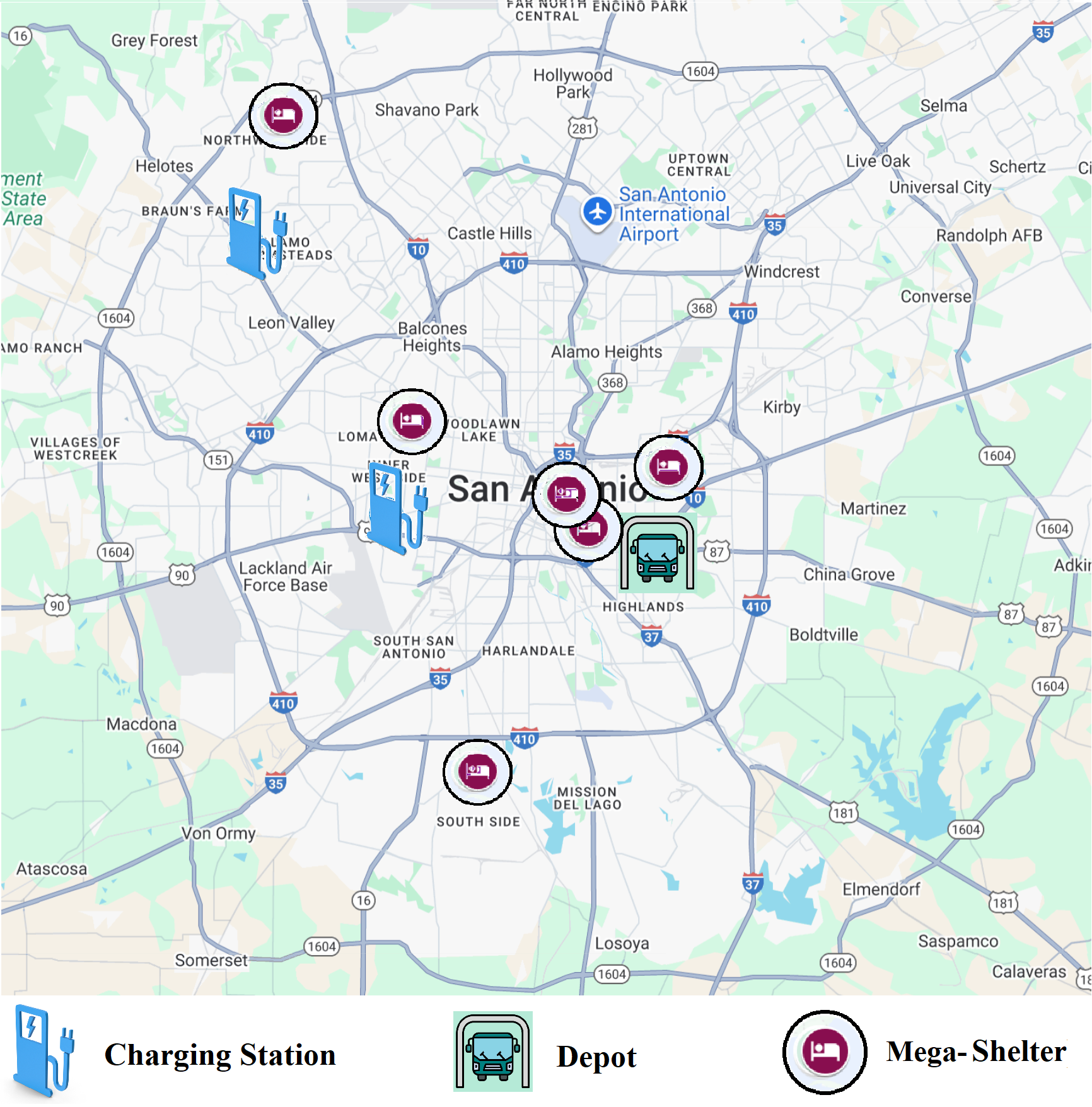}
\par\end{centering}
\caption{Geographical locations of mega-shelters, CSs, and the depot in the case study.\label{fig:Geographical_Locations}}
\end{figure}

In this section, we present a real case study to evaluate the performance of the proposed model and solution approaches in solving the \textit{ESBILR} for realistic-sized instances. Based on our close collaborations with the SAFD, we selected the six largest disaster shelters in San Antonio, Texas, United States, that are used by the SAFD as mega-shelters with large capacity and high energy demands during disasters. One of these shelters is used as a medical treatment facility for evacuees requiring medical treatment. 
We considered the San Antonio Independent School District Transportation as the depot for the ESB fleet in our case study. Due to the lack of available charging facilities for ESBs in San Antonio (there is no operational charging facility for ESBs yet), we considered two SAFD stations as CSs for the ESBs. Figure \ref{fig:Geographical_Locations} shows the geographical locations of the depot, mega-shelters, and CSs in our case study. Furthermore, to demonstrate the computational performance of our proposed models and solution algorithms in solving the \textit{ESBILR} for larger instances, we considered an additional SAFD station as a CS, as well as four additional critical facilities in San Antonio (that are water recycling centers and elementary schools) as critical isolated loads in our case study. We note that the schools are used as community-level shelters (with less capacity than mega-shelters) during disasters, and the water recycling centers are considered as critical facilities by SAFD.

\subsection{ESB Data}\label{subsec:ESB_Data}

ESBs are generally classified as long-/extended-range or fast-charge types, depending on the size of their battery packs. Long-/extended-range ESBs have battery capacities ranging from 250 to 660 kWh, whereas fast-charge ESBs' battery capacities vary between 50 to 250 kWh. 
A typical ESB, such as the Solobus, has a battery capacity of 300 kWh and an average speed of 40 mph (65 km/h)~\cite{gomotive_ebus_speed}. The electricity consumption of ESBs is strongly affected by the heating system,  
reported to span between 1 and 1.4 kWh/km on buses with fossil-fueled heating systems, and up to 2.35 kWh/km
on electrically heated ones~\cite{sustainablebus_ebus_energy}. 
The price of ESBs ranges from \$250,000 to \$440,000 based on the battery capacity~\cite{realclearenergy_esb_price}. Therefore, in this case study, we consider three different ESB types with the information shown in Table \ref{tab:ESB_Types_Information}. The energy consumption rate for the ESBs of type 3 is determined based on the average speed of 65 km/h and the average energy consumption of 1.18 kWh/km.

\begin{table}[h!]
\centering
\caption{Information of different ESB types.}
\label{tab:ESB_Types_Information}
\begin{tabular}{cccccc}
\hline
ESB Type & Cost (\$) & $C_{max}$ (kWh) & $C_{min}$ (kWh) & $G_{min}$ (kWh) & {\makecell[c]{Energy Consumption \\ Rate (kWh/h)}}\\ 
\hline \hline
Type 1 & 250,000 & 100 & 10 & 10 & 38.19 \\

Type 2 & 350,000 & 300 & 30 & 30 & 53.46 \\

Type 3 & 450,000 & 500 & 50 & 50 & 76.38 \\
\hline
\end{tabular}
\end{table}

\subsection{Service Times of Charging Stations and Shelters}\label{subsec:Time_Interval_Service_Times}

In this case study, we assume that the planning horizon comprises discrete time slots, each slot with a 15-minute duration, and thus we translate all the timings, including the service and travel times, to the number of time slots.
Fast charging can provide power up to 500 kW and is the quickest way to charge an EV, such as an ESB, in less than 30 minutes~\cite{edemands, he2020optimal}. Therefore, we consider the service time of CSs between 1 to 2 time slots in our case study. Regarding the service time of shelters, a majority of lithium-based batteries have a 1C discharging rate~\cite{lithium_discharge}, allowing a 300 kWh ESB to deliver up to 75 kWh per 15-minute time slot. Given the average demand of shelters between 100 to 250 kW per time slot, the service time of shelters should be assumed between 1 to 3 time slots. Therefore, we assume the service time of the shelters and CSs as shown in Table \ref{tab:Service_Time_Shelters_CSs}.

\begin{table}[h!]
\centering
\caption{Service time of shelters and CSs.}
\label{tab:Service_Time_Shelters_CSs}
\begin{tabular}{cccccccccccccc}
\hline
\multirow{2}{*}{Node}  & \multicolumn{10}{c}{Shelters} & \multicolumn{3}{c}{Charging Stations} \\
\cmidrule(lr){2-11} \cmidrule(lr){12-14}
& 1 & 2 & 3 & 4 & 5 & 6 & 7 & 8 & 9 & 10 & CS 1 & CS 2 & CS 3\\
\hline
\hline
Service Time (time slots)   & 1 & 1 & 2 & 3 & 2 & 3 & 1 & 1 & 3 & 3 & 1 & 2 & 1 \\
\hline
\end{tabular}
\end{table}

\subsection{Travel Times and Electricity Consumption}\label{subsec:Data_Description}
The detailed information of the travel times of ESBs between different nodes in the network is presented in Table \ref{tab:combined_travel_time} of Appendix \ref{subsec:Travel_Times}. The hourly demand profile of the shelters over a 48-time slots (i.e., 12-hour) planning horizon is demonstrated in Figure \ref{fig:Demand_Curve_Line} of Appendix \ref{subsec:Demand_Profile}.

\section{Numerical Results and Managerial Insights}\label{sec:Numerical_Results}

In this section, we assess the efficacy of the proposed models and algorithms, and compare the performance of different solution approaches in solving the \textit{ESBILR}. Moreover, based on the real case study presented in Section \ref{sec:Case_Study}, we conducted numerical experiments to provide practical insights for emergency management agencies into the following questions: (1) How do different solution approaches, Gurobi, B\&P with exact subproblem, and B\&P with labeling-algorithm-integrated dynamic programming heuristic subproblem, compare to each other in terms of runtime and solution quality (i.e., optimality gap)? (2) What effect does network sparsity level (i.e., shelter-ESB type compatibility) have on the performance of different solution approaches, the total cost, and the required ESB fleet composition? (3) What effect do the number of available ESBs, the average energy demand level of the shelters, the shelters' service time, and different levels of disaster severity have on the total cost, and the required ESB fleet composition and fleet size? and (4) How does the inconvenience fee of load shifting affect the percentage of shifted demand of the shelters, required fleet composition, and the total cost?

\subsection{Experimental Setup}\label{subsec:Experimental_Setup}

We implemented our models and solution methods in Python 3.12.7. We used the Gurobi solver ver. 12.0.1~\cite{gurobi2023} to solve the \textit{ESB-BaM-MILP} (i.e., Eqs. (\ref{eq:Obj Function}) - (\ref{eq:Continuous_Variables})), \textit{ESB-AccM-MILP} (i.e., Eqs. (\ref{eq:Obj Function}) - (\ref{eq:Continuous_Variables}) and (\ref{eq:Symmetry Breaking-1})), \textit{RLMP} (\ref{eq:master_problem}), model (\ref{eq:Subproblem}), and \textit{MILP-MP} (\ref{eq:MILP_master_problem}). We set the optimality gap and the time limit to 1\% and 1 hour, respectively, in all numerical experiments. We ran the numerical experiments on the compute nodes contained in the University of Texas at San Antonio High Performance Computing (HPC) Cluster, ARC, using a 64-bit Linux operating system. A compute node has two Intel Xeon Gold 6248 CPUs, each with @2.5 GHz and 40 physical cores with 376 GB of installed RAM. Based on parameter tuning, we solve the \textit{MILP-MP} (\ref{eq:MILP_master_problem}) after exploring the first 10 nodes of the B\&B tree and subsequently after every 20 nodes. For the heuristic B\&P approach, as it employs a heuristic for solving the subproblems, we terminate the B\&P if the lower bound improvement is less than 5\% over five successive nodes. In the model, we used $M_1 = 1000$, and we set $M_2$ to be two times larger than the total demand of all shelters over the planning horizon (e.g., $M_2=153,000$ for 10 shelters and a 48-time slots (i.e., 12-hour) planning horizon). 
For brevity, we categorize the setups used in our numerical experiments into distinct problem instances, denoted in the form S-C-T, representing the number of shelters, CSs, and time slots, respectively.

\subsection{Runtime and Solution Quality of Solution Approaches}\label{subsec:Runtime_Comparison}

In this section, Table \ref{tab:Solution_Approach_Comparison_Base_Case} shows the average runtime and optimality gap of the proposed models and solution approaches in solving the \textit{ESBILR} for different problem instances, obtained by running four replications of each problem instance with varying demand levels of the shelters. For brevity, we refer to solving \textit{ESB-BaM-MILP} (i.e., Eqs. (\ref{eq:Obj Function}) - (\ref{eq:Continuous_Variables})) and \textit{ESB-AccM-MILP} (i.e., Eqs. (\ref{eq:Obj Function}) - (\ref{eq:Continuous_Variables}) and (\ref{eq:Symmetry Breaking-1})) using Gurobi as \textit{BaM} and \textit{AcM}, respectively. The acronyms, \textit{BaBnP} and \textit{AcBnP}, denote solving the \textit{ESBILR} using the exact B\&P with and without solving the \textit{MILP-MP} (\ref{eq:MILP_master_problem}), respectively. We refer to the heuristic B\&P integrating dynamic programming and labeling algorithm to solve the subproblems as \textit{LaBnP}. The optimality gap of the \textit{LaBnP} is calculated as the percentage by which its objective value exceeds the optimal objective value of the \textit{AcBnP} approach.

\begin{table}[h!]
    \centering
    \caption{Runtime and solution quality comparison among different models and solution approaches. "-" and "*" represent that the corresponding solution approach is not able to compute a feasible solution within the one-hour time limit, and the solution quality for the corresponding problem instance cannot be reported, respectively. Problem instances are denoted as S-C-T, representing the number of shelters, CSs, and time slots (e.g., 1-1-16 represents one shelter, one CS, and 16 time slots).}
    \resizebox{\textwidth}{!}{
        \begin{tabular}{lccccccccccc}
            \toprule
            \multirow{3}{*}{\makecell[c]{Problem \\ Instance}} 
            & \multicolumn{5}{c}{Runtime (seconds)} 
            & \multicolumn{5}{c}{Optimality Gap (\%)} \\
            \cmidrule(lr){2-6} \cmidrule(lr){7-11}
            & \makecell[c]{\textit{BaM}}
            & \makecell[c]{\textit{AcM}}
            & \makecell[c]{\textit{BaBnP}}
            & \makecell[c]{\textit{AcBnP}}
            & \makecell[c]{\textit{LaBnP}}
            & \makecell[c]{\textit{BaM}}
            & \makecell[c]{\textit{AcM}}
            & \makecell[c]{\textit{BaBnP}}
            & \makecell[c]{\textit{AcBnP}}
            & \makecell[c]{\textit{LaBnP}} \\
            \midrule \hline
            1-1-16 
            & 1355 & 38 & 11 & 5 & 3 & 2.4 & 0.0 & 0.0 & 0.0 & 0.0 \\
            1-1-32 
            & 2861 & 211 & 62 & 10 & 6 & 14.5 & 0.0 & 0.0 & 0.0 & 0.0 \\
            1-1-48 
            & 2378 & 751 & 152 & 17 & 14 & 11.9 & 0.0 & 0.0 & 0.0 & 0.0 \\\hline
            2-1-16 
            & 3600 & 2538 & 199 & 22 & 5 & 31.6 & 0.0 & 0.0 & 0.0 & 0.0 \\
            2-1-32 
            & 3600 & 3600 & 807 & 131 & 14 & 49.3 & 36.7 & 0.0 & 0.0 & 0.0 \\
            2-1-48 
            & 3600 & 3600 & - & 310 & 35 & 61.8 & 58.8 & - & 0.0 & 1.6 \\\hline
            4-1-16 
            & 3600 & 3600 & - & 780 & 10 & 71.7 & 75.1 & - & 0.0 & 12.8 \\
            4-1-32 
            & 3600 & 3600 & - & 2302 & 40 & 88.6 & 87.5 & - & 0.0 & 11.9 \\
            4-1-48 
            & - & - & - & 3600 & 108 & - & - & - & 6.1 & 5.7 \\\hline
            6-2-16 
            & - & - & - & 2973 & 36 & - & - & - & 0.0 & 10.0 \\
            6-2-32 
            & - & - & - & - & 301 & - & - & - & - & * \\
            6-2-48 
            & - & - & - & - & 2371 & - & - & - & - & * \\\hline
            8-3-16 
            & - & - & - & 3600 & 128 & - & - & - & 5.7 & 9.9 \\
            8-3-32 
            & - & - & - & - & 2018 & - & - & - & - & * \\
            8-3-48 
            & - & - & - & - & - & - & - & - & - & - \\\hline
            10-3-16 
            & - & - & - & 3600 & 160 & - & - & - & 5.8 & 12.4 \\
            10-3-32 
            & - & - & - & - & 3554 & - & - & - & - & * \\
            10-3-48 
            & - & - & - & - & - & - & - & - & - & - \\
            \bottomrule
        \end{tabular}%
    }
    \label{tab:Solution_Approach_Comparison_Base_Case}
\end{table}

It is seen from Table \ref{tab:Solution_Approach_Comparison_Base_Case} that the performance of different solution approaches varies significantly across the problem instances. Among the Gurobi-based approaches (i.e., \textit{BaM} and \textit{AcM}), the symmetry-breaking constraints \eqref{eq:Symmetry Breaking-1} result in a better runtime performance of \textit{AcM} compared to \textit{BaM} for the small-scale problem instances. However, both approaches cannot solve the \textit{ESBILR} for the problem instances with more than two shelters and 32 time slots within the one-hour time limit. 

Among the exact B\&P-based approaches (i.e., \textit{BaBnP} and \textit{AcBnP}), integrating \textit{MILP-MP} (\ref{eq:MILP_master_problem}) in the B\&P procedure results in a better runtime performance of \textit{AcBnP} over \textit{BaBnP}. Moreover, \textit{AcBnP} solves a larger number of problem instances faster than the Gurobi-based approaches. On average, across all problem instances where \textit{AcM} converges to optimality within one hour, \textit{AcBnP} is \textit{47 times faster} than \textit{AcM}. Additionally, we observe from Table \ref{tab:Solution_Approach_Comparison_Base_Case} that, while the average runtime of \textit{BaM} for the small-scale problem instances comprising one shelter remains within one hour, its average optimality gap exceeds 1\%. This is because, the runtime of \textit{BaM} for some replications is much smaller compared to the one-hour time limit, whereas the \textit{BaM} cannot solve other replications of the same problem instance to optimality within the one-hour time limit. Although we varied the average demand parameter by very small amounts to create the four different replications of each problem instance, the runtime of the Gurobi-based approaches (e.g., \textit{BaM}) vary substantially across different replications of the same problem instance. On the other hand, the runtime of the B\&P-based approaches (e.g., \textit{AcBnP}) are consistent across different replications of the same problem instance. As Gurobi-based approaches are very sensitive to small parameter changes whereas the proposed B\&P algorithms are not, our proposed B\&P algorithms can serve as more reliable decision-support tools than the Gurobi-based approaches for leveraging ESBs in disaster restoration. 
To solve large-scale instances of the \textit{ESBILR} faster, we used the \textit{LaBnP} approach, which consistently outperforms the other solution approaches in terms of runtime while maintaining a good solution quality. On average, across all problem instances that \textit{AcBnP} converges to optimality within one hour, \textit{LaBnP} is \textit{27 times faster} than \textit{AcBnP}.

As explained above, \textit{AcM}, \textit{AcBnP}, and \textit{LaBnP} are the most effective solution approaches within their respective categories. To further evaluate their performance under varying network configurations, we define new instances of the \textit{ESBILR} for different shelter-ESB type compatibility using different network sparsity levels (SLs). We define four SLs from SL 1 (fully connected), where all ESB types can be plugged into all shelters, to SL 4 (fully sparse), where each shelter can only be connected to a single ESB type. 
The incidence matrices for different SLs are shown in Appendix \ref{Appendix_Incidence_Matrices}. The comparison of runtime and optimality gap among \textit{AcM}, \textit{AcBnP}, and \textit{LaBnP} in solving the \textit{ESBILR} for different SLs and problem instances is shown in Table \ref{tab:Solution_Approach_Comparison_Different_Sparisity_Levels_1} in Appendix \ref{Appendix_Sparsity_Results}.

We see from Table \ref{tab:Solution_Approach_Comparison_Different_Sparisity_Levels_1} that our proposed B\&P-based approaches (i.e., \textit{AcBnP} and \textit{LaBnP}) demonstrate a significantly better performance in terms of runtime over \textit{AcM} for different SLs. On average, across all SLs and problem instances where \textit{AcM} converges to optimality within one hour, \textit{AcBnP} and \textit{LaBnP} are, respectively, \textit{121 times} and \textit{335 times} faster than \textit{AcM}. Additionally, as the SL increases, the \textit{ESBILR} generally becomes more computationally challenging for the Gurobi-based approach (i.e., \textit{AcM}), whereas it becomes considerably easier for the B\&P-based approaches. This is because, at higher SLs, the \textit{ESBILR} solution requires a larger number of ESBs, increasing the number of decision variables for \textit{AcM} due to the explicit ESB indexing in Gurobi-based approaches, whereas the B\&P-based approaches generate ESB routes on an as-needed basis. Moreover, as the SL increases, each ESB type can serve fewer shelters, reducing the size of the CG subproblem and substantially decreasing the runtime of the B\&P-based approaches. Therefore, our B\&P algorithms better exploit the problem features and structure, enabling faster computations and handling larger problem instances across varying SLs, equipping emergency management agencies with faster decision-support tools to make reliable and accurate decisions in leveraging ESBs during disasters.

\subsection{Shelter-ESB Type Compatibility Vs. Total Cost and Fleet Composition}\label{subsec:Effect_Shelter_ESB_Matching}

\begin{figure}[H]
    \centering
    \begin{subfigure}[b]{0.48\textwidth}
        \centering
        \includegraphics[width=\textwidth]{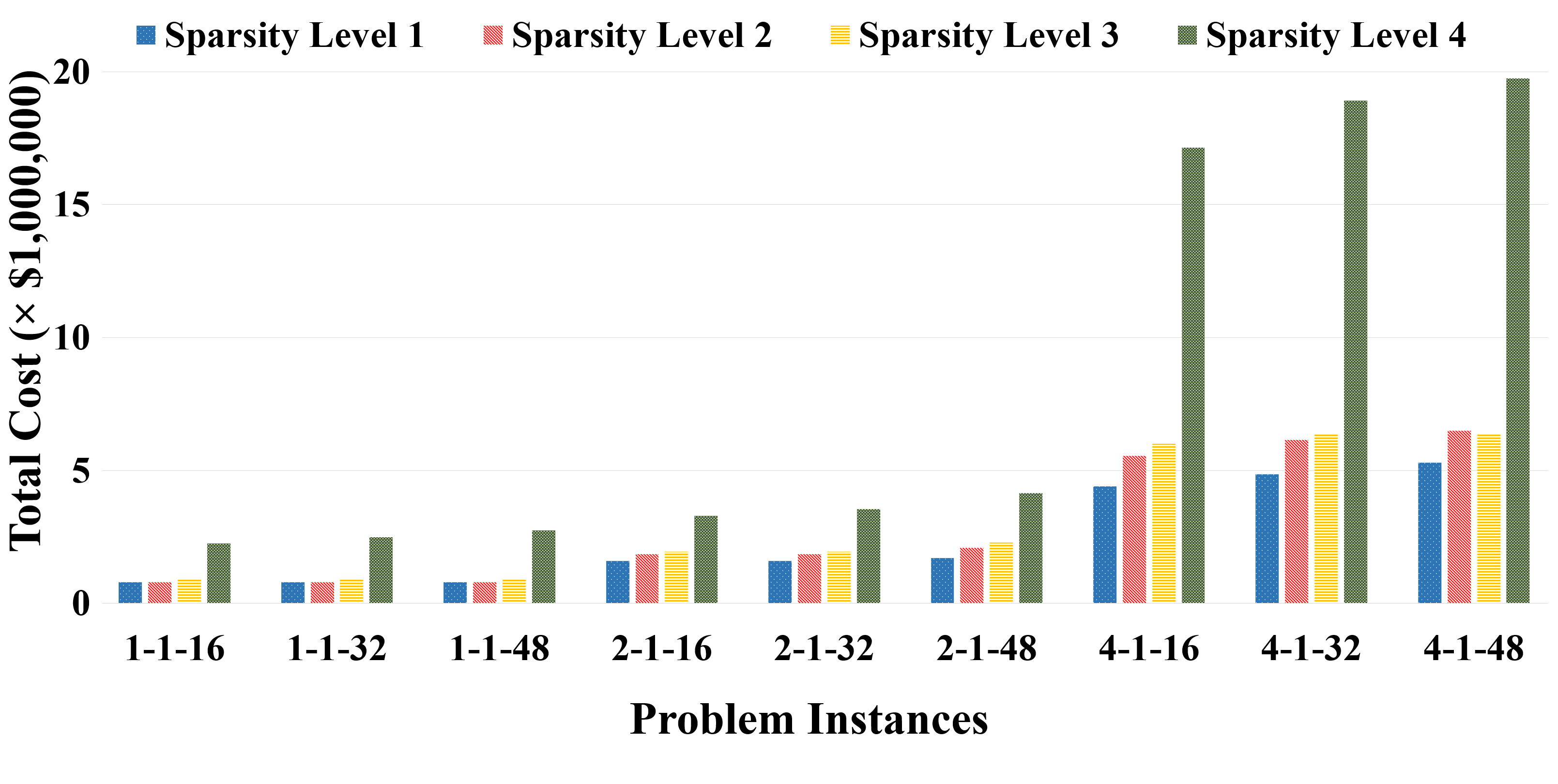}
        \caption{Effect of SL on the total cost.}
        \label{fig:Effect_Sparsity_on_Cost}
    \end{subfigure}
    \begin{subfigure}[b]{0.48\textwidth}
        \centering
        \includegraphics[width=\textwidth]{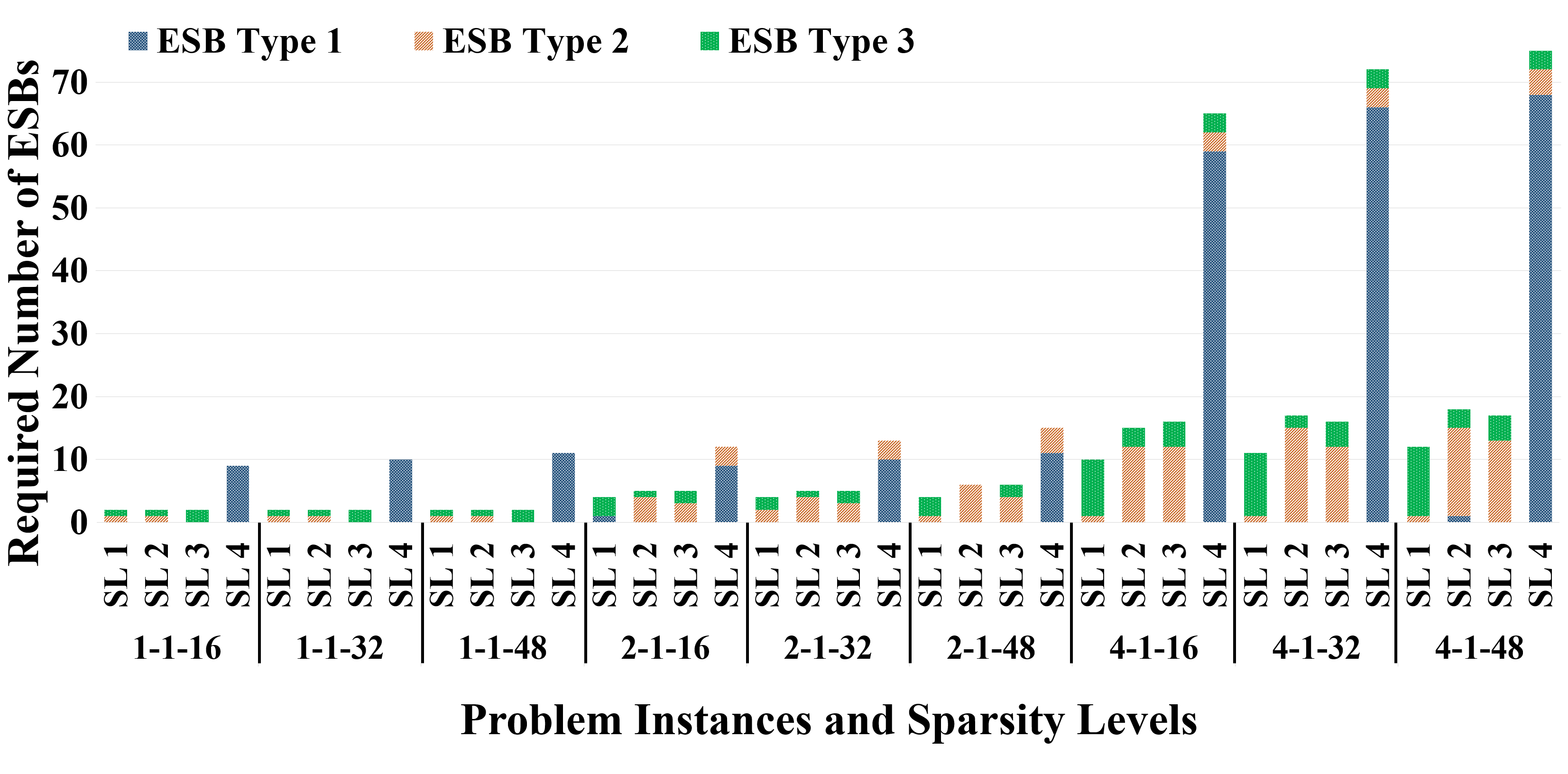}
        \caption{Effect of SL on the fleet composition.}
        \label{fig:Effect_Sparsity_on_ESB_Numbers}
    \end{subfigure}
    \caption{Effect of network sparsity level on the (a) total cost, and (b) required fleet composition.}
    \label{fig:Effect_Sparsity_Cost_Fleet}
\end{figure}

In this section, Figure \ref{fig:Effect_Sparsity_Cost_Fleet} demonstrates how different network sparsity levels (i.e., shelter-ESB type compatibility) arising from the technology, infrastructure, or spatial limitations affect the total cost, and the required fleet composition and fleet size in satisfying the energy demand of isolated loads.

We see from Figure \ref{fig:Effect_Sparsity_Cost_Fleet} that the total cost and the total number of required ESBs generally increase as the SL increases. At lower SLs, where the network is highly connected, each shelter can accommodate more than one ESB type, enabling the \textit{ESBILR} to leverage larger-capacity ESBs (i.e., type 2 and type 3) for all problem instances, as this is the most cost-efficient option. This reduces the total cost and the required number of ESBs as demonstrated in Figures \ref{fig:Effect_Sparsity_on_Cost} and \ref{fig:Effect_Sparsity_on_ESB_Numbers}, respectively. In contrast, higher SLs, where larger-capacity ESBs cannot be plugged into some shelters, force the optimal solution of \textit{ESBILR} to use a larger number of smaller-capacity ESBs (i.e., type 1) to satisfy the same energy demand, increasing the total cost and the required fleet size. Notably, on average across all the problem instances in this experiment, the total cost and the required fleet size increase by 207\% and 392\%, respectively, as the SL increases from level 1 to level 4.
This observation suggests that investing in shelter infrastructure to enhance the shelter-ESB type compatibility before disasters can reduce the total cost and the required fleet size, and thus improve overall efficiency. These findings are beneficial for emergency management agencies as they can benefit from these insights in designing and operating their ESB fleet to efficiently leverage ESBs for faster restoration of critical isolated loads around disasters.

\subsection{Number of Available ESBs Vs. Total Cost and Fleet Composition}\label{subsec:Effect_Available_ESBs}

Table \ref{tab:Effect_Power_Consumption} demonstrates the effect of the limited number of available ESBs on the total cost, and the required fleet composition and fleet size in satisfying the energy demand of isolated shelters. In this analysis, we gradually decrease the number of available ESBs, starting from no restriction on the number of available ESBs of types 2 and 3 (considered as the baseline) until the demand cannot be satisfied over the planning horizon with the available ESBs. 

\begin{table}[h!]
\centering
\caption{Effect of the number of available ESBs on the total cost and the required fleet composition of ESBs. "-" denotes that there is no limitation on the available number of ESBs of the corresponding type, whereas "*" denotes that the available ESBs cannot satisfy the energy demand of the shelters.}
\label{tab:Effect_Power_Consumption}
\begin{tabular}{ccccccc}
\hline
\multirow{2}{*}{\makecell[c]{Available \\ Type 2 ESBs}} & 
\multirow{2}{*}{\makecell[c]{Available \\ Type 3 ESBs}} & 
\multicolumn{2}{c}{Cost} & \multicolumn{3}{c}{Number of ESBs of Each Type} \\ \cmidrule(lr){3-4} \cmidrule(lr){5-7}
                                              &   & Total Cost (\$) & {\makecell[c]{Percentage \\ Increase}} & Type 1 & Type 2 & Type 3 \\ \hline \hline
 - & -   & 4,400,952          & 0\%           & 0    & 1    & 9  \\ 
- & 8    & 4,650,978          & 5.7\%         & 0    & 3    & 8  \\ 
- & 6    & 4,801,030          & 9.1\%        & 0    & 6    & 6  \\ 
- & 5    & 5,051,037          & 14.8\%        & 0    & 8    & 5  \\ 
- & 4    & 5,301,053          & 20.5\%        & 0    & 10    & 4  \\ 
- & 2    & 5,451,120          & 23.9\%        & 0    & 13   & 2  \\ 
- & 1    & 5,701,117          & 29.5\%        & 0    & 15   & 1  \\ 
- & 0    & 5,951,151          & 35.2\%        & 0    & 17   & 0  \\ 
16 & 0   & 6,351,222          & 44.3\%        & 3    & 16   & 0  \\
15 & 0   & 7,00,1347          & 59.1\%       & 7   & 15   & 0  \\ 
14 & 0   & 7,651,499          & 73.9\%       & 11   & 14   & 0  \\
13 & 0   & 8,301,656          & 88.6\%       & 15   & 13   & 0  \\ 
12 & 0   & 9,201,837          & 109.1\%       & 20   & 12   & 0  \\ 
8  & 0   & 12,302,507          & 179.5\%       & 38   & 8   & 0  \\ 
6 & 0   & *                & *             & *    & *    & *  \\ \hline
\end{tabular}
\end{table}

Table \ref{tab:Effect_Power_Consumption} shows that the total cost and the required fleet size increase as the availability of larger-capacity ESBs (i.e., type 2 and type 3 ESBs) decreases. 
When there are no restrictions on the availability of ESBs, the fleet predominantly comprises the most cost-efficient ESB type with the largest capacity (i.e., type 3 ESBs) to minimize the total cost. Decreasing the availability of larger-capacity ESBs reduces the flexibility of leveraging more cost-efficient ESBs with larger capacities, forcing the \textit{ESBILR} solution to use more ESBs of lower capacities to meet the energy demand of the shelters. Notably, it is seen from Table \ref{tab:Effect_Power_Consumption} that, on average, one unit reduction in the available number of the largest-capacity ESBs (i.e., type 3 ESBs) results in almost two unit increments in the required medium-capacity ESBs (i.e., type 2 ESBs). Moreover, reducing the availability of medium-capacity ESBs (i.e., type 2 ESBs) in a fleet that does not have any largest-capacity ESBs results in almost four unit increments in the required lowest-capacity ESBs (i.e., type 1 ESBs). As the cost of four additional type 1 ESBs is higher than the two additional type 2 ESBs, the rate of cost growth increases more significantly due to reducing the available number of medium-capacity ESBs compared to reducing the available number of the largest-capacity ESBs. For instance, we see from Table \ref{tab:Effect_Power_Consumption} that reducing the availability of type 3 ESBs by one unit from six to five results in a 5.2\% increase in the total cost, whereas reducing the availability of type 2 ESBs, in addition to no available type 3 ESBs, by one unit from 16 to 15 results in a 10.2\% increase in the total cost.

To quantify the number of required additional lower-capacity ESBs due to the reduction of a larger-capacity ESB in the fleet, we defined an \textbf{effective usable capacity} metric, $\zeta^k$, for ESB type $k$ computed using Eq. (\ref{eq:Effective_Usable_Capacity}). This metric represents how much energy, on average, an ESB of type $k$ has available in its battery to discharge to a shelter, accounting for its nominal capacity and the energy consumption during travel between different nodes in the network, ensuring the ESB can reach either a depot or a CS after serving the shelter.

\begin{equation}\label{eq:Effective_Usable_Capacity}
    \zeta^k = C_{max}^k - C_{min}^k - 2 \times T_{avg} \times \Gamma^{k}
\end{equation}

In Eq. (\ref{eq:Effective_Usable_Capacity}), $\Gamma^{k}$ is the energy consumption rate of ESB type $k$, and $T_{avg}$ is the average travel time between the depot and a shelter, as well as between a shelter and a CS. The factor 2 in $2 \times T_{avg} \times \Gamma^{k}$ accounts for the energy consumption of ESB type $k$ due to traveling from either the depot or a CS to a shelter, as well as returning to either the depot or a CS after serving a shelter. This ensures that the ESBs maintain a valid SOC at all times and have sufficient energy after serving the shelters to reach the depot or a CS. It is seen from Table \ref{tab:combined_travel_time} (in Appendix \ref{subsec:Travel_Times}) that $T_{avg} \simeq 0.47$ hours for the problem instance comprising four shelters and one CS used in this experiment. Therefore, according to Eq. (\ref{eq:Effective_Usable_Capacity}) and based on the available battery capacity (i.e., $C_{max} - C_{min}$) and the energy consumption rate (i.e., $\Gamma$) of the three ESB types shown in Table \ref{tab:ESB_Types_Information} of Section \ref{subsec:ESB_Data}, the \textit{effective usable capacity} of type 1, type 2, and type 3 ESBs are 54 kWh, 220 kWh, and 378 kWh, respectively. To reiterate, the largest-capacity ESB type has an average 378 kWh of available energy when arriving at a shelter, which can then be discharged to satisfy the energy demand of the shelter. In contrast, the lowest-capacity ESB type can only contribute to satisfying 54 kWh of the energy demand of shelters. The $\zeta^2/\zeta^1 = 220/54 \simeq 4$ reflects the required additional four ESBs of type 1 as the availability of type 2 decreases by one unit in addition to no availability of type 3 ESB. Similarly, $\zeta^3/\zeta^2 = 378/220\simeq 2$ reflects the required additional two ESBs of type 2 as the availability of type 3 decreases by one unit.          

This metric can help emergency management agencies obtain an estimation of the required number of ESBs of a given type to satisfy the energy demand of shelters. Also, it provides a relative measure of the number of additional ESBs of a lower-capacity ESB in case a larger-capacity ESB is unavailable. This metric is particularly useful for emergency management agencies seeking to determine the number of ESBs to dispatch for serving very few shelters for a short planning horizon when a comprehensive routing and scheduling model is not required.

To further highlight the significant difference in the cost-efficiency of the ESB types, we introduce a \textbf{capacity cost} metric, $\Upsilon^k$, for each ESB type $k$ computed using Eq. (\ref{eq:Capacity_Cost}). This metric represents the investment cost associated with each kWh of the \textit{effective usable capacity} for a given ESB type, which provides a measure of the cost-efficiency of an ESB type. The larger the value of this \textit{capacity cost} metric for an ESB type, the less cost-efficient that ESB type is.

\begin{equation}\label{eq:Capacity_Cost}
    \Upsilon^k = \frac{C_{inv}^k}{\zeta^k}
\end{equation}

According to Eq. (\ref{eq:Capacity_Cost}), the \textit{capacity cost} for type 1, type 2, and type 3 ESBs are \$4,630/kWh, \$1,591/kWh, and \$1,190/kWh, respectively. The computed values of the \textit{capacity cost} for the three ESB types used in our case study highlight the significant difference in the cost-efficiency of larger-capacity ESBs (i.e., type 2 and type 3 ESBs) compared to the lowest-capacity ESBs (i.e., type 1 ESBs). The significantly higher \textit{capacity cost} of type 1 ESBs better explains why type 1 ESBs---the least cost-efficient ESB type---are rarely used in the optimal fleet composition of most of the problem instances, whereas the optimal fleet composition predominantly comprises type 3 ESBs as the most cost-efficient ESB type. The \textit{capacity cost} serves as a metric to easily compare the relative benefit from different ESB types and evaluate their cost-effectiveness for meeting the energy demand of isolated shelters in disastrous situations. Emergency management agencies can use this metric to make better strategic investments in designing fleet composition under budgetary and operational limitations.

\subsection{Energy Demand Vs. the Total Cost and the Required Fleet Composition}\label{subsec:Effect_Average_Demand}

\begin{figure}[h!]
\begin{centering}
\includegraphics[scale=0.20]{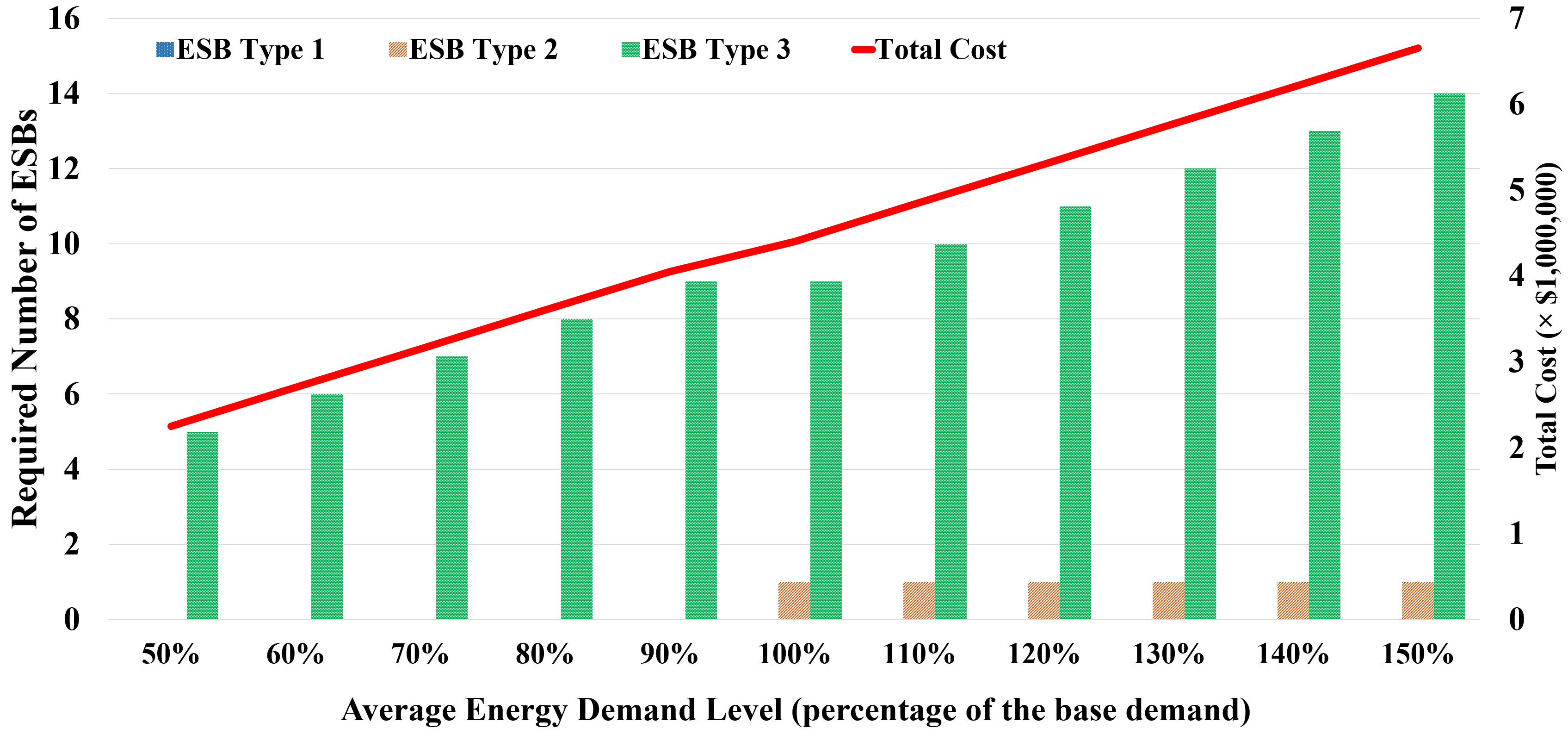}
\par\end{centering}
\caption{Effect of shelters' energy demand on total cost and required fleet composition. \label{fig:Effect_Avg_Demand_on_Cost}}
\end{figure}

This section presents the effect of increasing the energy demand of the shelters on the total cost and the required ESB fleet composition. Figure \ref{fig:Effect_Avg_Demand_on_Cost} illustrates a clear increasing trend in both the total cost and the required number of ESBs as the average energy demand increases from 50\% to 150\% of the base case demand. Specifically, the total cost, depicted by the red line, rises almost linearly with increasing demand, reflecting the growing need for additional resources to serve higher energy requirements. Notably, the required fleet size grows uniformly by one ESB, which results in the rate of cost increment of approximately \$440,000, for every 10\% increase in average demand (equal to an additional 1,636 kWh energy demand for the problem instance comprising four shelters, one CS, and 32 time slots used in this experiment).

Figure \ref{fig:Effect_Avg_Demand_on_Cost} also demonstrates the changes in the optimal fleet composition (number of ESBs of each type in the fleet) as the average energy demand increases. We see from Figure \ref{fig:Effect_Avg_Demand_on_Cost} that the fleet composition primarily comprises the largest-capacity ESBs (i.e., type 3 ESBs) as the average energy demand increases. This is due to the highest cost-efficiency of the largest-capacity ESBs (i.e., type 3 ESBs), as demonstrated in Section \ref{subsec:Effect_Available_ESBs} using the \textit{capacity cost} metric, compared to lower-capacity ESB types (i.e., types 1 and 2 ESBs). Notably, at the lower-level energy demands (e.g., 50\% to 90\%), the fleet comprises only the largest-capacity ESBs (i.e., type 3) for satisfying the energy demand of isolated shelters. Medium-capacity ESBs (i.e., type 2) are used less frequently and their contribution to the optimal fleet composition is usually much lower than the largest-capacity ESBs. In fact, the medium-capacity ESBs are used in the fleet along with the largest-capacity ESBs only when the requirement of additional capacity, given the current fleet composition, is not large enough to require a largest-capacity ESB incurring a high cost. In contrast to the largest-capacity ESBs, the lowest-capacity ESBs (i.e., type 1) are not used at any demand level due to their lowest cost-efficiency compared to the larger-capacity ESBs.

We defined a \textbf{capacity utilization} metric, $\varpi$, computed using Eq. (\ref{eq:Capacity_Utilization}), to highlight the benefit of our proposed \textit{ESBILR} problem and model in efficiently utilizing an ESB fleet with multiple back-and-forth trips between shelters and CSs to satisfy large energy demands of critical loads. This metric represents the number of times, on average, an ESB is visiting shelters during the planning horizon to satisfy a given amount of energy demand. The larger the value of the \textit{capacity utilization} metric, the higher the contribution of each ESB is in satisfying the demand, and therefore, the higher the efficiency of our proposed routing and scheduling model is. As computed in Section \ref{subsec:Effect_Available_ESBs}, the \textit{effective usable capacity} of type 3 and type 2 ESBs are 378 kWh and 220 kWh, respectively, with the cost of \$450,000 and \$350,000, respectively. As seen from Figure \ref{fig:Effect_Avg_Demand_on_Cost} and discussed in the above paragraph, every 10\% increase in the average energy demand (i.e., equal to 1,636 kWh) results in requiring an additional ESB in the fleet, more frequently of type 3 along with less frequent type 2, with a cost increase of approximately \$440,000. Using this information, we compute the \textit{effective usable capacity}, for an average ESB, $\bar{\zeta}$, that satisfies the 1,636 kWh increase in demand as Eq. (\ref{eq:Equivalent_Capacity}).

\begin{equation}\label{eq:Equivalent_Capacity}
    \begin{aligned}
        \text{average cost increment} &= \sum_{k \in \mathcal{K}} \varrho^k \times C_{inv}^k \\
        \bar{\zeta} &= \sum_{k \in \mathcal{K}} \varrho^k \times \zeta^k
    \end{aligned}
\end{equation}

\[
\$440,000 = 0.9 \times \$450,000 + 0.1 \times \$350,000 \quad \Rightarrow \quad 0.9 \times 378 + 0.1 \times 220 = 362 \, \text{kWh}
\]
\begin{equation}\label{eq:Capacity_Utilization}
\varpi = \frac{\text{energy demand}}{\bar{\zeta}}
\end{equation}

In Eq. (\ref{eq:Equivalent_Capacity}), $\varrho^k$ is the cost contribution of ESB of type $k$ in the average cost increment as the average energy demand increases such that $\sum_{k \in \mathcal{K}} \varrho^k = 1$. Therefore, the \textit{effective usable capacity} of an additional required ESB to satisfy the 1,636 kWh of the increased energy demand of the shelters is computed as $\bar{\zeta} = 362$ kWh. Using Eq. (\ref{eq:Capacity_Utilization}), the \textit{capacity utilization} of an ESB is $\varpi = 1,636/362 \simeq 4.5$, which means, on average, an ESB performs between four and five back-and-forth trips during the planning horizon, effectively satisfying 4.5 times higher energy demand than its \textit{effective usable capacity}. We observed that this \textit{capacity utilization} is consistently higher than one across different problem instances in our real case study. This higher \textit{capacity utilization} is due to the efficient ESB routing and scheduling solution provided by the proposed model and algorithms, where each ESB performs multiple back-and-forth trips between shelters and CSs, which allows each ESB to deliver more energy than its \textit{effective usable capacity} through multiple trips.     

\subsection{Shelters' Service Time Vs. Total Cost and Fleet Composition}\label{subsec:Effect_Service_Time}

In this section, Table \ref{tab:Effect_Service_Time} demonstrates the effect of increasing the service time of the shelters on the total cost and the required ESB fleet composition. In this analysis, we consider the duration of each time slot to be 5 minutes. We increase the service time from one time slot (considered as the baseline) to six time slots to demonstrate the effect of the shelters' service time on the total cost and the required ESB fleet composition. We see an increasing trend (in Table \ref{tab:Effect_Service_Time}) in the total cost and the required fleet size as the service time increases. Notably, as the service time increases from the shortest service time (i.e., 5 minutes) to the longest service time (i.e., 30 minutes), the optimal solution of \textit{ESBILR} requires an additional type 3 ESB---the most cost-efficient ESB type with the least \textit{capacity cost}---resulting in a 36\% rise in the total cost. This is because, the ESBs discharge energy more quickly with shorter service times, enabling more frequent back-and-forth trips between shelters and CSs, and therefore, making better utilization of the ESBs in the critical load restoration. In contrast, longer service times reduce the frequency of back-and-forth trips of ESBs, reducing the \textit{capacity utilization} of each ESB and requiring a larger fleet size (i.e., more ESBs) with higher capacities to satisfy the same total energy demand of shelters.

\begin{table}[h!]
\centering
\caption{Effect of service time of the shelters on the total cost and the required fleet composition.}
\label{tab:Effect_Service_Time}
\begin{tabular}{cccccc}
\hline
\multirow{2}{*}{Service Time (minutes)} & \multicolumn{2}{c}{Cost} & \multicolumn{3}{c}{Number of ESBs of Each Type} \\ \cmidrule(lr){2-3} \cmidrule(lr){4-6}
                                        & Total Cost (\$) & Percentage Increase & Type 1 & Type 2 & Type 3 \\ \hline \hline
5                                       & 1,250,322         & -           & 0      & 1      & 2      \\ 
15                                      & 1,350,316         & 8.0\%       & 0      & 0      & 3      \\ 
20                                      & 1,600,302         & 28.0\%      & 0      & 2      & 2      \\ 
30                                      & 1,700,319         & 36.0\%      & 0      & 1      & 3      \\ \hline
\end{tabular}
\end{table}

\subsection{Disaster Severity Vs. Total Cost and Fleet Composition}\label{subsec:Effect_Disaster_Severity}

In this section, Table \ref{tab:Effect_Disaster_Severity} shows the effect of increasing the severity level of the disaster on the total cost and the required fleet composition. In this analysis, we defined three disaster severity levels by changing the travel times and energy consumption rate of the ESBs (to account for the effect of weather severity and traffic congestion) as follows. (1) Normal weather - base case travel time and energy consumption rate, (2) moderate weather - travel times increased by one time slot and energy consumption rates increased by 20\%, and (3) adverse weather - travel times are doubled and the energy consumption rates increased by 50\%.

\begin{table}[ht]
\centering
\caption{Effect of disaster severity on the total cost and the required fleet composition of ESBs.}
\label{tab:Effect_Disaster_Severity}
\begin{tabular}{cccccc}
\hline
\multirow{2}{*}{Disaster Severity} & \multicolumn{2}{c}{Cost}                     & \multicolumn{3}{c}{Number of ESBs of Each Type} \\ \cmidrule(lr){2-3} \cmidrule(lr){4-6} 
                                    & Total Cost (\$)           & Percentage Increase       & Type 1 & Type 2 & Type 3 \\ \hline \hline
Normal Weather                      & 4,400,952            & -                         & 0      & 1      & 9      \\ 
Moderate Weather                    & 6,652,022            & 51.1\%                    & 0      & 1      & 14     \\ 
Adverse Weather                     & 14,306,529           & 225.0\%                   & 0      & 1      & 31     \\ \hline
\end{tabular}
\end{table}

We see from Table \ref{tab:Effect_Disaster_Severity} that, the total cost and the required fleet size significantly increase as the severity level of the disaster increases. Notably, increasing the severity level of the disaster from normal weather to adverse weather results in a 225\% increase in the total cost and more than three times larger required ESB fleet comprising the largest capacity ESB type (i.e., type 3 ESBs). This increment in the total cost and the fleet size is due to longer travel times and higher energy consumption rates of ESBs caused by the higher severity level of the disaster. Specifically, longer travel times reduce the frequency of back-and-forth trips of ESBs, reducing the \textit{capacity utilization} of ESBs (i.e., the number of times each ESB can visit shelters within the planning horizon to discharge energy), as discussed in Section \ref{subsec:Effect_Average_Demand}. Moreover, according to Eq. (\ref{eq:Effective_Usable_Capacity}) in Section \ref{subsec:Effect_Available_ESBs}, longer travel times and higher energy consumption rates reduce the \textit{effective usable capacity} of ESBs, reducing the amount of energy each ESB has available in its battery to discharge to a shelter. Therefore, the compounded effect of the reduced \textit{capacity utilization} and \textit{effective usable capacity} of ESBs due to the increased disaster severity level forces the \textit{ESBILR} solution to use a larger fleet size (i.e., more ESBs) comprising larger-capacity ESBs to satisfy the same total energy demand of isolated shelters, resulting in a higher total cost.

The effect of disaster severity on the required fleet size and total cost does not increase uniformly as the disaster severity level increases; rather, it increases exponentially. This is because of the compounded effect of longer travel times and higher energy consumption rates on the required fleet size and the total cost as explained above. For instance, the number of required additional ESBs is 5 as the weather changes from normal to moderate, whereas the number of required additional ESBs is 17 (i.e., more than 3 times higher than additional ESBs from normal to moderate) as the weather changes from moderate to adverse.

\subsection{Inconvenience Fee of Demand Shifting Vs. Percentage of Shifted Demand, Total Cost, and Fleet Composition}\label{subsec:Effect_Inconvenience_Fee}

\begin{figure}[H]
    \centering
    \begin{subfigure}[b]{0.48\textwidth}
        \centering
        \includegraphics[width=\textwidth]{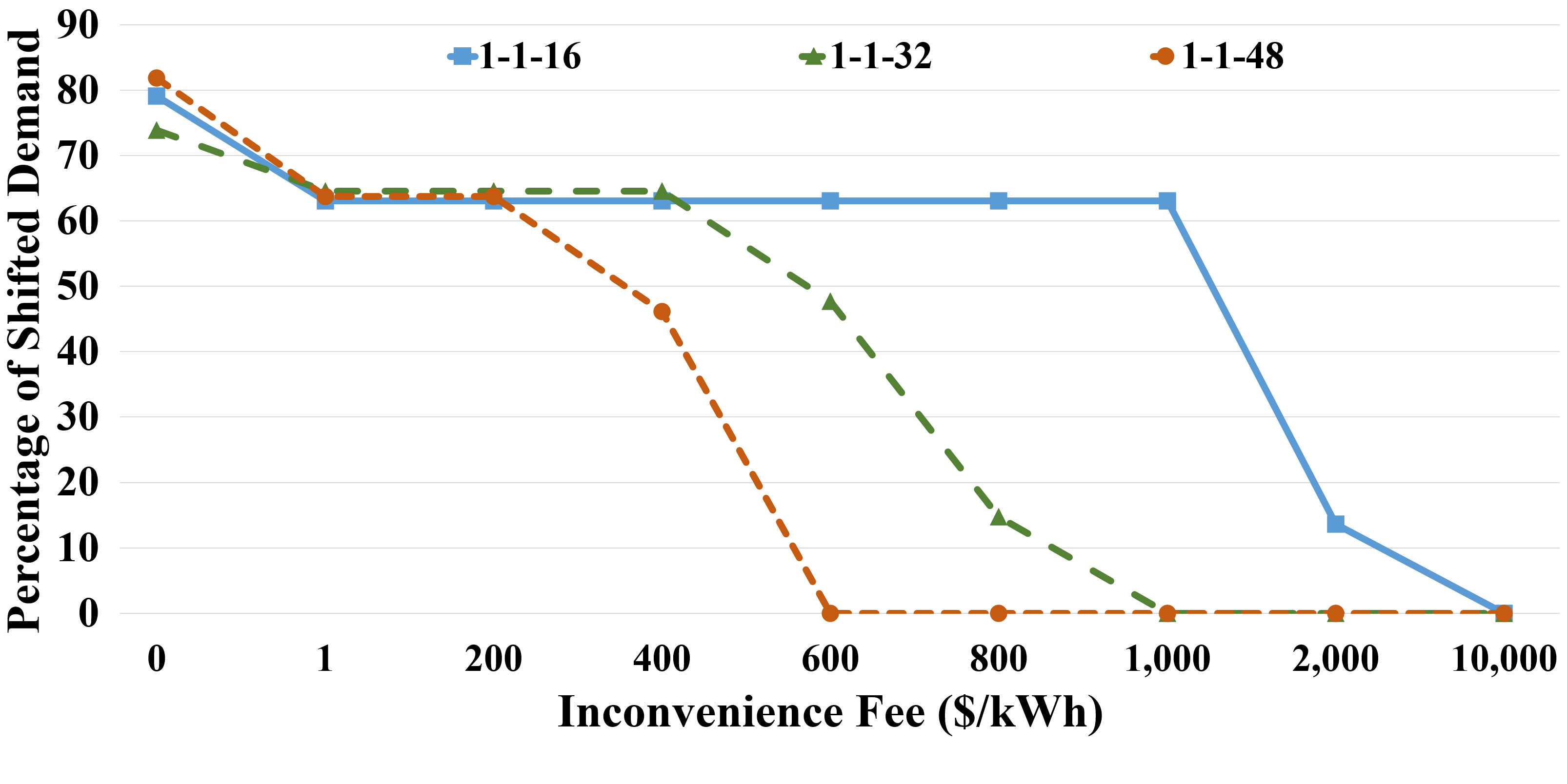}
        \caption{Percentage of shifted demand.}
        \label{fig:Effect_inconv_on_shifted_demand}
    \end{subfigure}
    \begin{subfigure}[b]{0.48\textwidth}
        \centering
        \includegraphics[width=\textwidth]{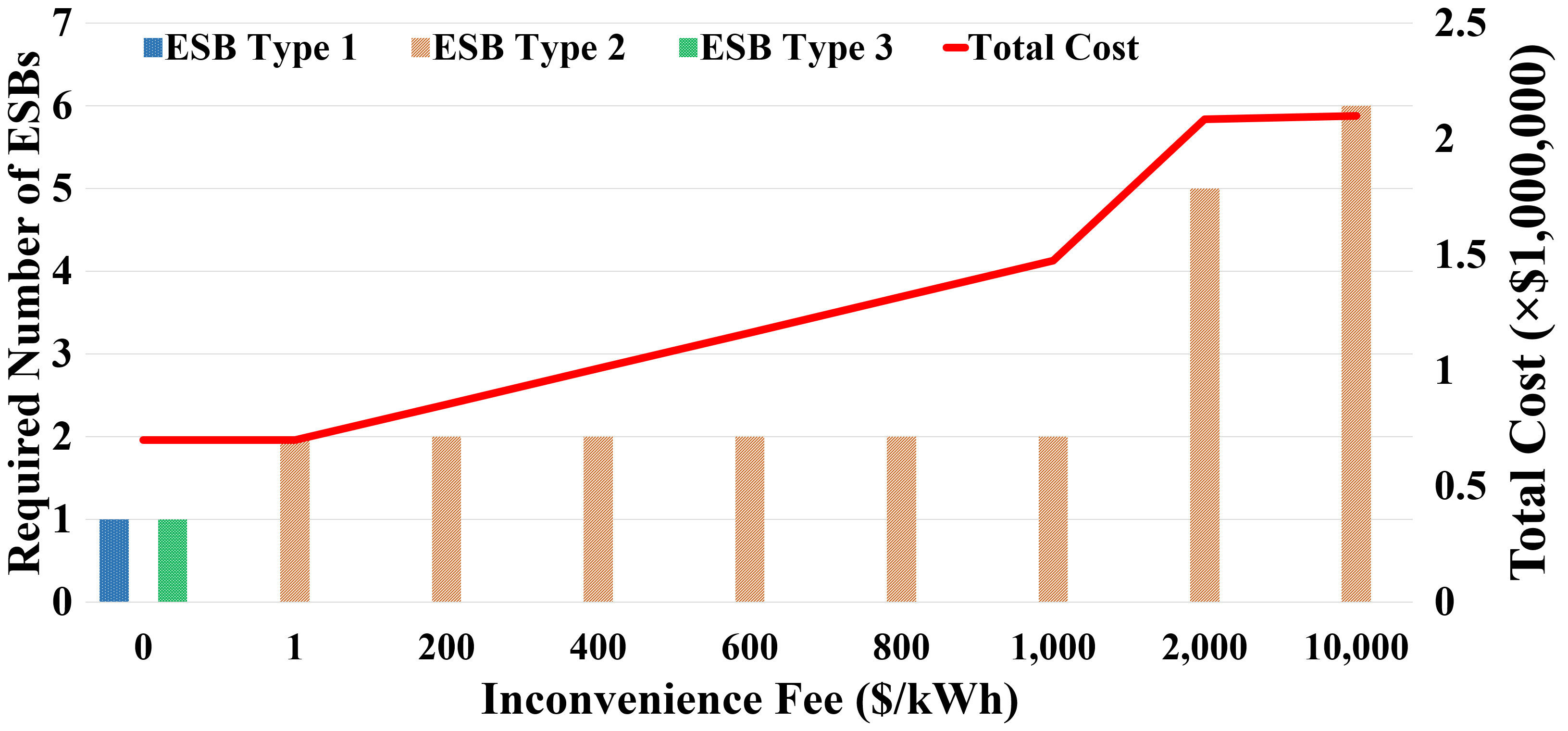}
        \caption{Total cost and fleet composition for 1-1-16.}
        \label{fig:Effect_inconv_on_Cost_ESB}
    \end{subfigure}
    \caption{Inconvenience fee vs. (a) shifted demand, and (b) total cost and fleet composition.}
    \label{fig:Effect_Inconvenience_Demand_Shifting_Total_Cost}
\end{figure}

In this section, we model the reluctance of the shelters to shift their demand from one time slot to others using the inconvenience fee (cost) of demand shifting. Figures \ref{fig:Effect_inconv_on_shifted_demand} and \ref{fig:Effect_inconv_on_Cost_ESB} show the effect of increasing the inconvenience fee of demand shifting on the percentage of shifted demand, and the total cost and the required fleet composition, respectively. To model the inconvenience of shifting energy demand (i.e., load) from one time slot to the other time slots, we define $l_{i,t}^{'}$ and $F_{i,t}^{'}$ as the amount of shifted demand and the inconvenience fee of shifting the demand of shelter $i$ at time slot $t$, respectively. We modify the \textit{ESB-AccM-MILP} (i.e., Eqs. (\ref{eq:Obj Function}) - (\ref{eq:Continuous_Variables}) and (\ref{eq:Symmetry Breaking-1})), \textit{RLMP} (\ref{eq:master_problem}), and model (\ref{eq:Subproblem}) as explained in Appendix \ref{Appendix_Inconvenience_Modeling}. As no standardized value exists for inconvenience fees in disastrous situations and selecting an appropriate value for the inconvenience fee depends on various factors (e.g., the type of the shelter), in this analysis, we vary the inconvenience fee from no inconvenience fee (i.e., $F_{i,t}^{'}=\$0/kWh$), representing complete flexibility in shifting the demand to an extremely high inconvenience fee (i.e., $F_{i,t}^{'}=\$10,000/kWh$), representing a strict situation where no demand can be shifted to other time slots.

We see from Figure \ref{fig:Effect_inconv_on_shifted_demand} that, as the inconvenience fee increases, the reluctance of the shelter to shift the energy demand across the planning horizon increases, decreasing the percentage of shifted demand. 
Figure \ref{fig:Effect_inconv_on_Cost_ESB} shows that the total restoration cost and the required fleet size increase as the inconvenience fee increases for all problem instances (refer to Figure \ref{fig:Effect_inconv_on_Cost_ESB_2_3} in Appendix \ref{Appendix_Inconvenience_Extra_Results} for problem instances 1-1-32 and 1-1-48). Notably, as the inconvenience fee increases to extremely high values, the total cost and the required fleet size increase significantly by up to 200\%.
This is because, at extremely high inconvenience fees, the contribution of the inconvenience cost of demand shifting in the total restoration cost, defined in Eq. (\ref{eq:Obj_Function_Shifting_MILP}) in Appendix \ref{Appendix_Inconvenience_Modeling}, outweighs the contribution of the investment cost of the ESB fleet in the total cost (in the objective function). Therefore, the optimal solution of \textit{ESBILR} uses more ESBs in the fleet to strictly satisfy the time-dependent energy demand (i.e., the percentage of shifted demand reaches zero), increasing the total cost and the required fleet size.
These findings provide valuable insights for practitioners into how the inconvenience fee 
affects the total cost and required fleet composition. Practitioners can use the proposed models and algorithms on their own data, ESB fleet, and an appropriate inconvenience fee based on their specific needs
to determine the required ESB fleet composition and fleet size, as well as their routes and schedules.

\section{Conclusion}\label{sec:Conclusion}

In this paper, we study the problem of routing and energy scheduling of a heterogeneous fleet of government-owned ESBs with V2B technology within a discretized planning horizon for faster power restoration to critical isolated loads around disasters while minimizing the total investment and transportation energy consumption costs. We address several practical aspects of the problem, including simultaneous transportation and energy scheduling of ESBs, route interdependency and spatial-wise coupling among ESBs, multiple back-and-forth trips between shelters and CSs with partial energy delivery. We modeled the problem as an efficient MILP formulation enhanced with valid inequalities. We proposed customized exact and labeling-algorithm-integrated dynamic programming heuristic B\&P algorithms to efficiently solve the power restoration problem. Numerical results based on a real case study of disaster shelters in San Antonio, Texas, United States, demonstrate that our proposed exact B\&P algorithm is \textit{121 times faster} than the Gurobi solver and the heuristic B\&P provides high-quality solutions while computationally \textit{27 times faster} than the exact B\&P algorithm. 

Using network sparsity, we modeled the limited service capability of ESBs due to technology, infrastructure, or spatial limitations. Results highlight that our proposed B\&P-based approaches better exploit the network sparsity feature than Gurobi.
The \textit{effective usable capacity} and \textit{capacity cost} metrics indicate that the largest-capacity ESBs are seven times more effective and four times more cost-efficient than the least-capacity ESBs. The \textit{capacity utilization} metric shows that each ESB satisfies 4.5 \textit{times higher} demand than its \textit{effective usable capacity}, highlighting the efficiency of our proposed ESB routing and scheduling approach.
Results also show that as the weather changes from normal to adverse, the total cost increases by 225\%. Furthermore, results demonstrate that as the inconvenience fee of demand shifting increases from no inconvenience fee (\$0/kWh) to an extremely high inconvenience fee (\$10,000/kWh), the percentage of shifted demand reduces to zero, whereas the total cost and fleet size increase by up to 200\%.

The outcomes of this research are pivotal for the practical implementation of ESBs in faster power restoration to isolated shelters during disasters. Our proposed B\&P algorithms can serve as fast and reliable decision-support tools for emergency management agencies in leveraging ESBs around disasters. The proposed \textit{capacity cost} and \textit{effective usable capacity} metrics help emergency management agencies in assessing the cost-efficiency of different ESB types and the required number of ESBs of a given type for satisfying the energy demand of isolated shelters under varying weather conditions. While the efficacy of the optimization model and the algorithms were demonstrated with specific ESB types, problem settings, and the novel practical aspects, practitioners can use the proposed algorithms on their own data and network configuration to design and operate their ESB fleet in faster power restoration to critical isolated loads around disasters.

This study can be extended to account for the uncertainties in the demand profile of the shelters, weather conditions, and the availability of roads in a stochastic or robust optimization framework. Another future direction is to study $N-h$ (i.e., ESB outage during the restoration) and $N+i$ (i.e., increment in the number of shelters) contingencies as a two-stage robust optimization model, where the first-stage model determines the required ESB fleet composition, whereas the routing and energy scheduling of ESBs are determined in the second stage under the worst-case contingencies.

\bibliographystyle{unsrt}
\bibliography{References_Full}

\appendix 

\newpage
\setcounter{page}{1}
\renewcommand{\thefigure}{\Alph{section}.\arabic{figure}}
\renewcommand{\thetable}{\Alph{section}.\arabic{table}}
\renewcommand{\theequation}{\Alph{section}.\arabic{equation}}
\renewcommand{\thepage}{A-\arabic{page}}
\renewcommand{\thealgorithm}{\Alph{section}.\arabic{algorithm}}

\newpage

\section*{Appendices}

\setcounter{equation}{0}
\setcounter{figure}{0}
\setcounter{table}{0}
\section{Proofs}
\label{Appendix_Proof}

\renewcommand{\theprop}{\ref{prop:symmetry_breaking}}

\begin{prop}
    The \textit{ESB-BaM-MILP} model augmented with symmetry-breaking constraints (\ref{eq:Symmetry Breaking}) - (\ref{eq:Symmetry Breaking-1}) remains as a valid model, and its feasible set is a proper subset of that of \textit{ESB-BaM-MILP}.
\end{prop}

\begin{proof}
    Let $F$ be the set of feasible solutions of the original problem without symmetry-breaking constraints (i.e., \textit{ESB-BaM-MILP} model defined in Eqs. (\ref{eq:Obj Function}) - (\ref{eq:Continuous_Variables})) and $F'$ as the set of feasible solutions after adding symmetry-breaking constraints (\ref{eq:Symmetry Breaking}) and (\ref{eq:Symmetry Breaking-1}). As multiple ESBs of the same type exist, there can be multiple optimal solutions that differ only in the ESBs used from the same ESB type. Let $h_1, h_2 \in \mathcal{H}^k$ be two ESBs of type $k$ such that $h_1 < h_2$, and suppose there exist two optimal solutions $x_1, x_2 \in F$ where (a) in $x_1$, ESB $h_1$ follows a route and energy schedule while ESB $h_2$ is not used, and (b) in $x_2$, ESB $h_2$ follows the same route and energy schedule while $h_1$ is not used. As the two ESBs are from the same ESB type and have the same routing and energy scheduling patterns, we conclude that $f(x_1) = f(x_2)$; i.e., both solutions result in the same objective value, creating redundant symmetry in the solution space. By enforcing the symmetry-breaking constraints (\ref{eq:Symmetry Breaking}) - (\ref{eq:Symmetry Breaking-1}), any solution $x_2$ where a higher-index ESB is assigned but a lower-index ESB is not, becomes infeasible. The remaining feasible solutions are only those where ESBs are activated in increasing index order, making $F' \subseteq F$. However, as the symmetry-breaking constraints (\ref{eq:Symmetry Breaking}) - (\ref{eq:Symmetry Breaking-1}) do not eliminate distinct optimal solutions that involve different routing, energy scheduling, and the required ESB fleet composition but only remove the symmetric equivalents (i.e., differ only in the index of ESBs), the optimal solution and objective value remain unchanged. Therefore, symmetry-breaking constraints (\ref{eq:Symmetry Breaking}) - (\ref{eq:Symmetry Breaking-1}) are valid inequalities that result in a reduction in the search space without affecting the optimal solution, improving the computational efficiency of the \textit{ESB-BaM-MILP}.  
\end{proof}

\renewcommand{\theprop}{\ref{prop:NP-hardness}}

\begin{prop}
    The \textit{ESBILR} is NP-hard.
\end{prop}

\begin{proof}
    Consider the Euclidean Traveling Salesman Problem (TSP) where a vehicle visits each node $v\in \mathcal{V}$ in a complete graph $G = (\mathcal{V}, \mathcal{E})$ defined with Euclidean distance while minimizing the total distance traveled. On top of $G$, we construct a special case of the \textit{ESBILR} problem as follows. We arbitrarily select one node $v^0\in \mathcal V$ as the depot, consider a single time slot, assign one unit of energy demand for other nodes, and define that the energy consumption of each edge is proportional to its length (say $kl_e$ with $l_e$ being the length of edge $e$ and $k$ being the proportional parameter). Moreover, consider a single ESB with a capacity greater than $|\mathcal V|+k\sum_{e\in \mathcal E}l_e$. 
    
    Clearly, once the ESB is fully charged at $v^0$, no additional charge is needed for its travel over all edges in $G$. Also, we assume it can discharge and move very fast so that those operations/movements can be completed in one time slot. Indeed, because of the triangle inequality and the cost minimization objective function, we can easily show that $(i)$ all nodes will be visited only once for discharge in any optimal solution of this \textit{ESBILR}; $(ii)$ an optimal route with the minimum cost is a route of the minimum distance.  
    Therefore, we can conclude that if \textit{ESBILR} is solved exactly, its optimal solution yields an optimal solution to the TSP.    
    
    Given that the TSP is NP-hard, it follows that \textit{ESBILR} is NP-hard too.    
\end{proof}

\renewcommand{\theprop}{\ref{prop_dominance}}

\begin{prop}
    Let $\mathbb A$ denote the set of feasible labels obtained so far, and consider two labels $L$ and $L'$ in $\mathbb A$.
    $(i)$ If $o(L)=o(L')$, $rc(L')<rc(L)$,  $soc(L')\geq soc(L),$ and $f(L')<f(L)$, we can discard $\{L\}$ from $\mathbb A$ without losing any optimal route. $(ii)$ If $L$ and $L'$ visit the same set of nodes, $rc(L')<rc(L),$ and $f(L')<f(L)$, we can discard $\{L\}$ from $\mathbb A$ without losing any optimal route.
\end{prop}

\begin{proof}
    The first statement can be easily proven by contradiction. Let $R(L)$ denote the route associated with label $L$. Also, assume that there exists an optimal full route obtained by extending label $L$ to $\hat L$. Hence, $R(\hat L)\backslash R(L)$ denotes a partial route starting from $o(L)$ to $l'$. Note that $R(L')\cup R(\hat L)\backslash R(L)$ is a new route between $l$ and $l'$. Actually, according to the listed assumptions in $(i)$ and Eqs. \eqref{L1} - \eqref{L5}, it not only is feasible but also has a smaller reduced cost, which contradicts to the assumption that $R(\hat L)$ is optimal. Hence, the desired conclusion follows. Regarding the second statement, we have $soc(L)=soc(L')$ given that their partial routes visit the same set of nodes (with different orders). The conclusion follows directly from the first statement.     
\end{proof}

\setcounter{equation}{0}
\setcounter{figure}{0}
\setcounter{table}{0}
\section{Detail Description of Case Study}
\label{Appendix_Case_Study_Details}

\subsection{Travel Times Between Different Locations}\label{subsec:Travel_Times}

We assume that all ESB types have the same travel time. 
The rationale is that, the discretization of the planning horizon makes small variations in the speed of ESBs negligible. 
Additionally, in a practical disaster scenario, the speeds of vehicles are usually limited by the road conditions (e.g., weather and congestion due to evacuation traffic), where different ESBs will take approximately the same time to travel between two nodes in the network, despite small differences in their nominal speed. The approximate travel times between different nodes in the network, calculated using the actual road distances and the average ESB speed of 65 km/h, are shown in Table \ref{tab:combined_travel_time}.

\begin{table}[htbp]
\centering
\caption{Travel time between shelters, depot, and charging stations (number of time slots).}
\label{tab:combined_travel_time}
\begin{tabular}{ccccccccccccccc}
\hline
\diagbox[dir=NW]{From}{To} & \multicolumn{10}{c}{Shelters} & \multicolumn{1}{c}{\multirow{2}{*}{Depot}} & \multicolumn{3}{c}{Charging Stations} \\ 
\cline{2-11} \cline{13-15}
 & 1 & 2 & 3 & 4 & 5 & 6 & 7 & 8 & 9 & 10 & & CS 1 & CS 2 & CS 3\\ 
\hline \hline
Shelter 1 & 0 & 1 & 1 & 3 & 2 & 2 & 3 & 2 & 2 & 3 & 1 & 2 & 3 & 2   \\
Shelter 2 & 1 & 0 & 1 & 3 & 2 & 2 & 3 & 3 & 3 & 3 & 1 & 2 & 3 & 2   \\
Shelter 3 & 1 & 1 & 0 & 3 & 2 & 2 & 3 & 2 & 2 & 3 & 1 & 2 & 2 & 2   \\
Shelter 4 & 3 & 3 & 3 & 0 & 2 & 3 & 5 & 3 & 2 & 2 & 3 & 3 & 1 & 2   \\
Shelter 5 & 2 & 2 & 2 & 2 & 0 & 2 & 4 & 2 & 2 & 2 & 2 & 1 & 2 & 1   \\
Shelter 6 & 2 & 2 & 2 & 3 & 2 & 0 & 2 & 2 & 3 & 3 & 2 & 2 & 3 & 2   \\
Shelter 7 & 3 & 3 & 3 & 5 & 4 & 2 & 0 & 3 & 4 & 4 & 3 & 3 & 4 & 4   \\
Shelter 8 & 2 & 3 & 2 & 3 & 2 & 2 & 3 & 0 & 2 & 2 & 2 & 2 & 2 & 2   \\
Shelter 9 & 2 & 3 & 2 & 2 & 2 & 3 & 4 & 2 & 0 & 2 & 3 & 3 & 2 & 1   \\
Shelter 10 & 3 & 3 & 3 & 2 & 2 & 3 & 4 & 2 & 2 & 0 & 3 & 2 & 1 & 2   \\
\hline
\end{tabular}
\end{table}

\subsection{Electricity Consumption and Demand Profile of Shelters}\label{subsec:Demand_Profile}

\begin{figure}[H]
\begin{centering}
\includegraphics[scale=0.233]{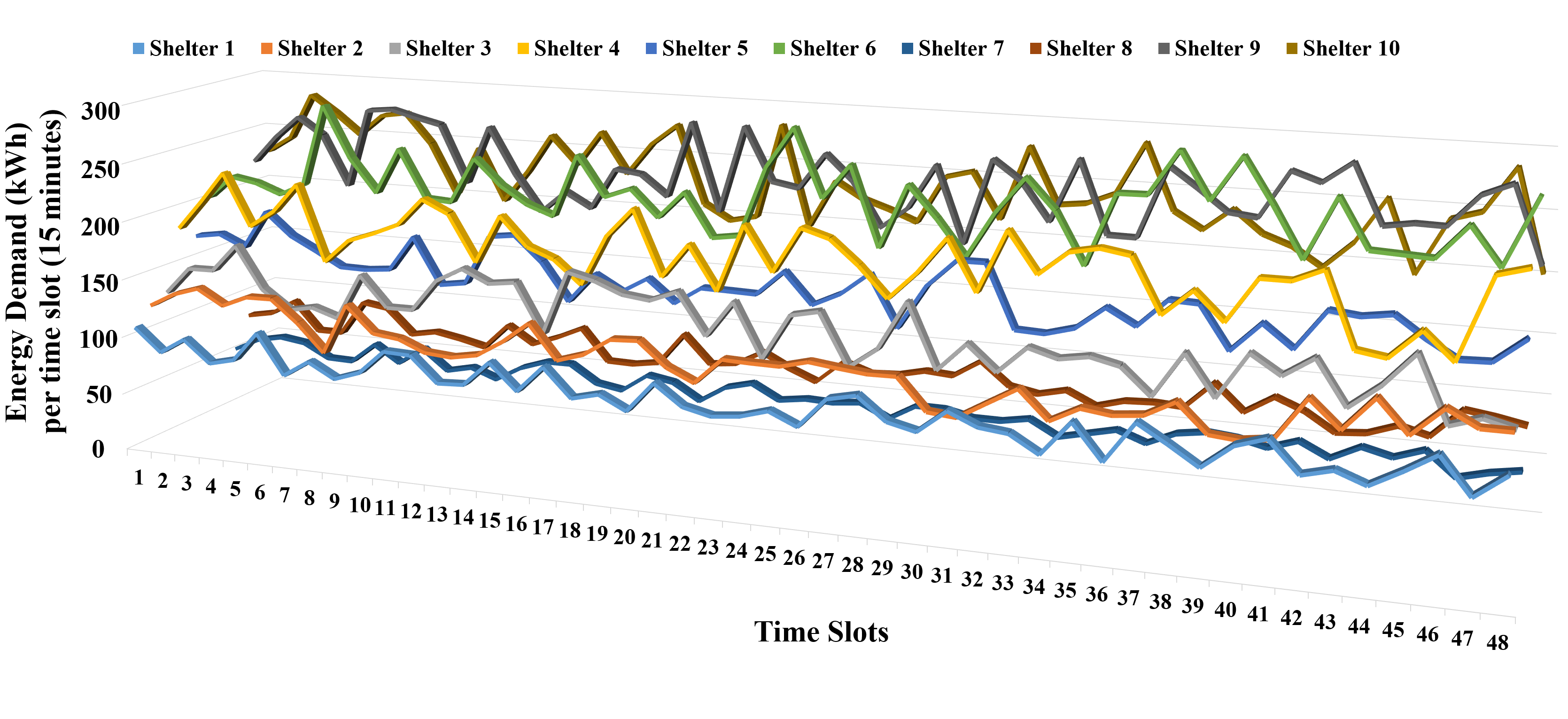}
\par\end{centering}
\caption{Demand profile of shelters over the planning horizon.\label{fig:Demand_Curve_Line}}
\end{figure}

According to our close collaborations with the SAFD, each mega-shelter is designed to accommodate more than 5,000 people. As each mega-shelter requires $40 \sim 60$ square feet per person \cite{Megashelter_guide}, the average area of each mega-shelter is around 250,000 square feet in the U.S. Furthermore, mega-shelters provide a range of resources for the people from basic needs to medical assistance, which requires $24 \sim 31$ kWh (i.e., 27.5 kWh on average) energy per square foot annually~\cite{electricityplans_hospitals, dsoelectric_hospitals, eia_commercial_2012, koupaei2023identifying}. Therefore, the estimated average hourly electricity consumption of a mega-shelter is $27.5 \times 250,000 / 8760 \approx 785\,\text{kWh}$. Due to the variability in the size of mega-shelters and their varied energy consumption over time, we assume the energy requirement of the mega-shelters falls between 100 to 250 kWh per 15-minute time slot. Regarding the water recycling centers, it is indicated that the recycled water system in San Antonio consumes an average of approximately 8.4 MWh per year, resulting in an average of 953 kWh energy consumption of water recycling facilities per hour (i.e., approximately 240 kWh per 15-minute time slot). Moreover, the two elementary schools in our case study have the area of 107,340 and 75,900 square feet, resulting in the estimated hourly average electricity consumption of 340 and 240 kWh, respectively (i.e., 85 and 60 kWh per 15-minute time slot, respectively). Figure \ref{fig:Demand_Curve_Line} presents the demand profile of the shelters over a 48-time slots (i.e., 12-hour) planning horizon. As the longest travel time from the depot to the shelters is three time slots, we set the demand of the initial three time slots to zero, ensuring all ESBs can start servicing the shelters from the beginning of the demand profile.

\setcounter{equation}{0}
\setcounter{figure}{0}
\setcounter{table}{0}
\section{Incidence Matrices for Different Sparsity Levels}
\label{Appendix_Incidence_Matrices}

We define four different SLs to model the limited shelter-ESB type compatibility. The SLs are modeled using incidence matrices (IM), which consist of binary elements, indicating the compatibility between shelters and ESB types. Each element $IM(k,i)$ is 1, if shelter $i$ is compatible with ESB type $k$ and 0, otherwise.
The incidence matrices for the four SLs are defined as follows.

\begin{align*}
\text{IM - SL 1} &=
\begin{bmatrix}
1 & 1 & 1 & 1 & 1 & 1 & 1 & 1 & 1 & 1 \\
1 & 1 & 1 & 1 & 1 & 1 & 1 & 1 & 1 & 1 \\
1 & 1 & 1 & 1 & 1 & 1 & 1 & 1 & 1 & 1 \\
\end{bmatrix},
&
\text{IM - SL 2} &=
\begin{bmatrix}
1 & 1 & 0 & 1 & 0 & 1 & 1 & 1 & 0 & 1 \\
1 & 1 & 1 & 1 & 1 & 0 & 1 & 1 & 1 & 0 \\
1 & 0 & 1 & 0 & 1 & 1 & 1 & 0 & 1 & 1 \\
\end{bmatrix}
\end{align*}

\begin{align*}
\text{IM - SL 3} &=
\begin{bmatrix}
1 & 1 & 0 & 1 & 0 & 1 & 1 & 0 & 0 & 1 \\
0 & 1 & 0 & 1 & 1 & 0 & 0 & 1 & 0 & 0 \\
1 & 0 & 1 & 0 & 0 & 1 & 1 & 0 & 1 & 1 \\
\end{bmatrix},
&
\text{IM - SL 4} &=
\begin{bmatrix}
1 & 0 & 0 & 1 & 0 & 0 & 1 & 0 & 0 & 1 \\
0 & 1 & 0 & 0 & 1 & 0 & 0 & 1 & 0 & 0 \\
0 & 0 & 1 & 0 & 0 & 1 & 0 & 0 & 1 & 0 \\
\end{bmatrix}
\end{align*}

\setcounter{equation}{0}
\setcounter{figure}{0}
\setcounter{table}{0}
\section{Performance Comparison Among Algorithms for Varying Sparsity Levels}
\label{Appendix_Sparsity_Results}

The numerical results of the runtime and optimality gap of \textit{AcM}, \textit{AcBnP}, and \textit{LaBnP} in solving the \textit{ESBILR} for different SLs and problem instances are shown in Table \ref{tab:Solution_Approach_Comparison_Different_Sparisity_Levels_1}. 

\begin{sidewaystable}[htbp]
    \centering
    \small
    \caption{Runtime and solution quality of different solution approaches for varying sparsity levels.}
    \setlength{\tabcolsep}{2pt}
    \begin{tabular}{cccccccccccccccc}
        \hline
        Sparsity & Problem & \multicolumn{3}{c}{Runtime (seconds)} & \multicolumn{3}{c}{Optimality Gap (\%)} 
        & Sparsity & Problem & \multicolumn{3}{c}{Runtime (seconds)} & \multicolumn{3}{c}{Optimality Gap (\%)} \\
        \cmidrule(lr){3-5} \cmidrule(lr){6-8} \cmidrule(lr){11-13} \cmidrule(lr){14-16}
        Level & Instance & \textit{AcM} & \textit{AcBnP} & \textit{LaBnP} & \textit{AcM} & \textit{AcBnP} & \textit{LaBnP} 
        & Level & Instance & \textit{AcM} & \textit{AcBnP} & \textit{LaBnP} & \textit{AcM} & \textit{AcBnP} & \textit{LaBnP} \\
        \cmidrule(lr){1-8} \cmidrule(lr){9-16} \addlinespace[1pt] 
        \cmidrule(lr){1-8} \cmidrule(lr){9-16}
        \multirow{18}{*}{1} & 1-1-16  & 14 & 4 & 3 & 0.0 & 0.0 & 0.0  & \multirow{18}{*}{3} & 1-1-16 & 89 & 2 & 1 & 0.0 & 0.0 & 0.0 \\
        & 1-1-32 & 99 & 10 & 6 & 0.0 & 0.0 & 0.0  &        & 1-1-32 & 1067 & 3 & 2 & 0.0 & 0.0 & 0.0 \\
        & 1-1-48 & 208 & 18 & 14 & 0.0 & 0.0 & 0.0  &     & 1-1-48 & 3600 & 13 & 4 & 4.6 & 0.0 & 0.0  \\
        & 2-1-16  & 2360 & 29 & 5 & 0.0 & 0.0 & 0.0  &   & 2-1-16 & - & 7 & 4 & - & 0.0 & 0.0  \\
        & 2-1-32 & 3600 & 175 & 14 & 52.2 & 0.0 & 0.0  &  & 2-1-32 & - & 31 & 10 & - & 0.0 & 0.0  \\
        & 2-1-48 & 3600 & 245 & 33 & 61.3 & 0.0 & 0.0  &  & 2-1-48 & - & 154 & 24 & - & 0.0 & 0.0  \\
        & 4-1-16  & 3600 & 184 & 10 & 71.1 & 0.0 & 11.4 &  & 4-1-16 & - & 23 & 6 & - & 0.0 & 0.0  \\
        & 4-1-32 & 3600 & 2265 & 41 & 86.9 & 0.0 & 9.3 &  & 4-1-32 & - & 170 & 23 & - & 0.0 & 7.9  \\
        & 4-1-48 & -    & 3600 & 112& -    & 7.2 & 1.9  &  & 4-1-48 & - & 715 & 58 & - & 0.0 & 7.1  \\
        & 6-2-16  & -    & 2093 & 36 & -    & 0.0 & 8.9 &  & 6-2-16 & - & 76 & 18 & - & 0.0 & 8.1  \\
        & 6-2-32 & -    & -   & 312& -    & -   & *    &  & 6-2-32 & - & 3512 & 153 & - & 0.0 & 4.1 \\
        & 6-2-48 & -    & -   & 2791&-    & -   & *    &  & 6-2-48 & - & 3600 & 956 & - & 3.4 & 3.7 \\
        & 8-3-16 & -    & 3600 & 117& -    & 4.6 & 8.6  &  & 8-3-16 & - & 3443 & 37 & - & 0.0 & 6.2 \\
        & 8-3-32& -    & -   & 1829& -   & -   & *    &  & 8-3-32 & - & 3600 & 772 & - & 8.9 & 5.3 \\
        & 8-3-48& -    & -   & -  & -    & -   & -    &  & 8-3-48 & - & - & 3600 & - & - & * \\
        & 10-3-16 & -    & 3600 & 152& -    & 5.6 & 8.3  &  & 10-3-16 & - & 3600 & 57 & - & 4.1 & 9.4 \\
        & 10-3-32& -    & -   & 3417& -   & -   & *    &  & 10-3-32 & - & 3600 & 1564 & - & 18.3 & 6.5 \\
        & 10-3-48& -    & -   & -  & -    & -   & -    &  & 10-3-48 & - & - & 3600 & - & - & * \\
        \cmidrule(lr){1-8} \cmidrule(lr){9-16}
        \multirow{18}{*}{2} & 1-1-16  & 58  & 4   & 3   & 0.0 & 0.0 & 0.0  & \multirow{18}{*}{4} & 1-1-16 & 133 & 2 & 1 & 0.0 & 0.0 & 0.0 \\
        & 1-1-32 & 674 & 9   & 6  & 0.0  & 0.0 & 0.0  &  & 1-1-32 & 2082 & 4 & 1 & 0.0 & 0.0 & 0.0 \\
        & 1-1-48 & 2536 & 17  & 12  & 0.0  & 0.0 & 0.0  &  & 1-1-48 & - & 11 & 1 & - & 0.0 & 0.0  \\
        & 2-1-16  & 3600 & 18  & 4   & 68.5 & 0.0 & 0.0  &  & 2-1-16 & - & 3 & 1 & - & 0.0 & 0.0  \\
        & 2-1-32 & -    & 111 & 13  & -    & 0.0 & 0.0  &  & 2-1-32 & - & 7 & 2 & - & 0.0 & 0.0  \\
        & 2-1-48 & -    & 405& 5  & -    & 0.0 & 4.8  &  & 2-1-48 & - & 18 & 5 & - & 0.0 & 0.0  \\
        & 4-1-16  & -    & 128 & 10  & -    & 0.0 & 6.3  &  & 4-1-16 & - & 7 & 3 & - & 0.0 & 0.0  \\
        & 4-1-32 & -    & 3563& 28  & -    & 0.0 & 11.4  &  & 4-1-32 & - & 43 & 13 & - & 0.0 & 0.0  \\
        & 4-1-48 & -    & 3600& 78 & -    & 4.0 & 3.1  &  & 4-1-48 & - & 322 & 35 & - & 0.0 & 0.0  \\
        & 6-2-16  & -    & 1618& 22  & -    & 0.0 & 5.7  &  & 6-2-16 & - & 35 & 10 & - & 0.0 & 4.4  \\
        & 6-2-32 & -    & 3600& 230 & -    & 6.2 & 0.0 &  & 6-2-32 & - & 401 & 98 & - & 0.0 & 4.5 \\
        & 6-2-48 & -    & -   & 1556& -    & -   & *    &  & 6-2-48 & - & 3359 & 725 & - & 0.0 & 2.9 \\
        & 8-3-16 & -    & 3600& 76  & -    & 1.3 & 11.0  &  & 8-3-16 & - & 160 & 23 & - & 0.0 & 7.2 \\
        & 8-3-32& -    & -   & 1765& -    & -   & *    &  & 8-3-32 & - & 795 & 572 & - & 0.0 & 2.7 \\
        & 8-3-48& -    & -   & 3600   & -    & -   & *    &  & 8-3-48 & - & 3600 & 3600 & - & 2.8 & 6.9 \\
        & 10-3-16 & -    & 3600& 86  & -    & 4.7 & 5.9  &  & 10-3-16 & - & 281 & 37 & - & 0.0 & 0.6 \\
        & 10-3-32& -    & -   & 2991& -    & -   & *    &  & 10-3-32 & - & 2233 & 785 & - & 0.0 & 7.2 \\
        & 10-3-48& -    & -   & -   & -    & -   & -    &  & 10-3-48 & - & 3600 & 3600 & - & 11.6 & 4.3 \\
        \hline
    \end{tabular}
    \label{tab:Solution_Approach_Comparison_Different_Sparisity_Levels_1}
\end{sidewaystable}

\setcounter{equation}{0}
\setcounter{figure}{0}
\setcounter{table}{0}
\section{Modeling Inconvenience Fee of Demand Shifting}
\label{Appendix_Inconvenience_Modeling}

To model the inconvenience of shifting energy demand from one time slot to the other time slots, we define $l_{i,t}^{'}$ and $F_{i,t}^{'}$ as the amount of shifted demand and the inconvenience fee of shifting the demand of shelter $i$ at time slot $t$, respectively. Then, we modify the objective function of the \textit{ESB-AccM-MILP} model (i.e., Eq. (\ref{eq:Obj Function})) as Eq. (\ref{eq:Obj_Function_Shifting_MILP}), and add constraint (\ref{eq:Load_Shifting_MILP_Model}) to the \textit{ESB-AccM-MILP} model.
\begin{align}
& min \qquad \Bigg\{\sum_{k\in\mathcal{K}}\sum_{h_k\in\mathcal{H}^k}\sum_{t\in\mathcal{T}}\sum_{i\in\mathcal{I}^k}\left(C_{inv}^{k}u_{0,i,t}^{h_k}\right)+ \sum_{i \in \mathcal{I}} F_i l_i + \sum_{i\in\mathcal{I}}\sum_{t\in\mathcal{T}}F_{i,t}^{'}l_{i,t}^{'}+\label{eq:Obj_Function_Shifting_MILP}\\
& \qquad \qquad
C_{fx}\sum_{k\in\mathcal{K}}\sum_{h_k\in\mathcal{H}^k}\sum_{i\in\mathcal{I}^k}\sum_{t\in\mathcal{T}}\left(\sum_{q\in\mathcal{Q}^k}\left(x_{i,q,t}^{h_k}R_{i,q}^{k}\right) + \sum_{j\in\mathcal{J}}\left(z_{j,i,t}^{h_k}R_{j,i}^{k}\right) + \left(u_{0,i,t}^{h_k}R_{0,i}^{k}\right) + \left(r_{i,0,t}^{h_k}R_{i,0}^{k}\right)\right)\Bigg\}\nonumber\\
& s.t. \qquad l_{i, t}^{'} \geq P_{i,t} - \sum_{k \in \mathcal{K}^i}\sum_{h_k\in\mathcal{H}^k}g_{i,t}^{h_k},\qquad\forall i\in\mathcal{I}, \forall t \in \mathcal{T} \label{eq:Load_Shifting_MILP_Model}
\end{align}
To incorporate the inconvenience of demand shifting in the B\&P approach, we made the following modifications: (1) the objective function (i.e., Eq. (\ref{eq:MPa})) of the master problem is changed to Eq. (\ref{eq:Obj_Function_Shifting_BnP}), (2) the constraint (\ref{eq:Load_Shifting_BnP_Model}) is added to the master problem with the dual variable $\delta_{i,t}$, and (3) the objective function of the subproblem (i.e., model (\ref{eq:Subproblem})) is changed to Eq. (\ref{eq:Subproblem_Objective_Shifting}).
\begin{align}
& min \qquad \sum_{k\in\mathcal{K}}\sum_{p_{k}\in\mathcal{P}_{k}}C_{p_{k}}^{k}\lambda_{p_{k}}^{k} + \sum_{i\in\mathcal{I}}F_{i}l_{i} + \sum_{i\in\mathcal{I}}\sum_{t\in\mathcal{T}}F_{i,t}^{'}l_{i,t}^{'} \label{eq:Obj_Function_Shifting_BnP}\\
& s.t. \qquad l_{i,t}^{'} \geq P_{i,t} - \sum_{k\in\mathcal{K}^i}\sum_{p_{k}\in\mathcal{P}_{k}}G_{i,t,p_{k}}\lambda_{p_{k}}^{k},\qquad\forall i\in\mathcal{I}, \forall t \in \mathcal{T}\label{eq:Load_Shifting_BnP_Model}
\end{align}

\begin{align}
\label{eq:Subproblem_Objective_Shifting}
& min \qquad \Bigg\{ C_{fx}\sum_{i\in\mathcal{I}^k}\sum_{t\in\mathcal{T}}\left(\sum_{q\in\mathcal{Q}}\left(x_{i,q,t}R_{i,q}^{k}\right) + \sum_{j\in\mathcal{J}}\left(z_{j,i,t}R_{j,i}^{k}\right) + \left(u_{0,i,t}R_{0,i}^{k}\right) + \left(r_{i,0,t}R_{i,0}^{k}\right)\right)\\\nonumber
& \qquad \qquad + \sum_{t\in\mathcal{T}}\sum_{i\in\mathcal{I}^k}C_{inv}^{k}u_{0,i,t}+\sum_{i\in\mathcal{I}^k}\sum_{t\in\mathcal{T}}w_{i,t}\left\{ \mu_{i}+\mu_{i}^{k}+\rho_{i,t}+\rho_{i,t}^{k}\right\} +\psi^{k}\\ \nonumber
& \qquad  +\sum_{i\in\mathcal{I}^k}\sum_{j\in\mathcal{J}}\left(x_{i,j,t}+z_{j,i,t}\right)\left\{ \eta_{i,j}+\eta_{i,j}^{k}\right\} +\sum_{i\in\mathcal{I}^k}\left(u_{0,i,t}+r_{i,0,t}\right)\left\{ \theta_{i}+\theta_{i}^{k}\right\} -\sum_{i\in\mathcal{I}}\sum_{t\in\mathcal{T}}\left(\left(\pi_{i}+ \delta_{i, t}\right)g_{i,t}\right) \Bigg\}\nonumber
\end{align}

\setcounter{equation}{0}
\setcounter{figure}{0}
\setcounter{table}{0}
\section{Inconvenience Fee Vs. Total Cost and Fleet Composition }
\label{Appendix_Inconvenience_Extra_Results}

\begin{figure}[H]
    \centering
    \begin{subfigure}[b]{0.48\textwidth}
        \centering
        \includegraphics[scale=0.11]{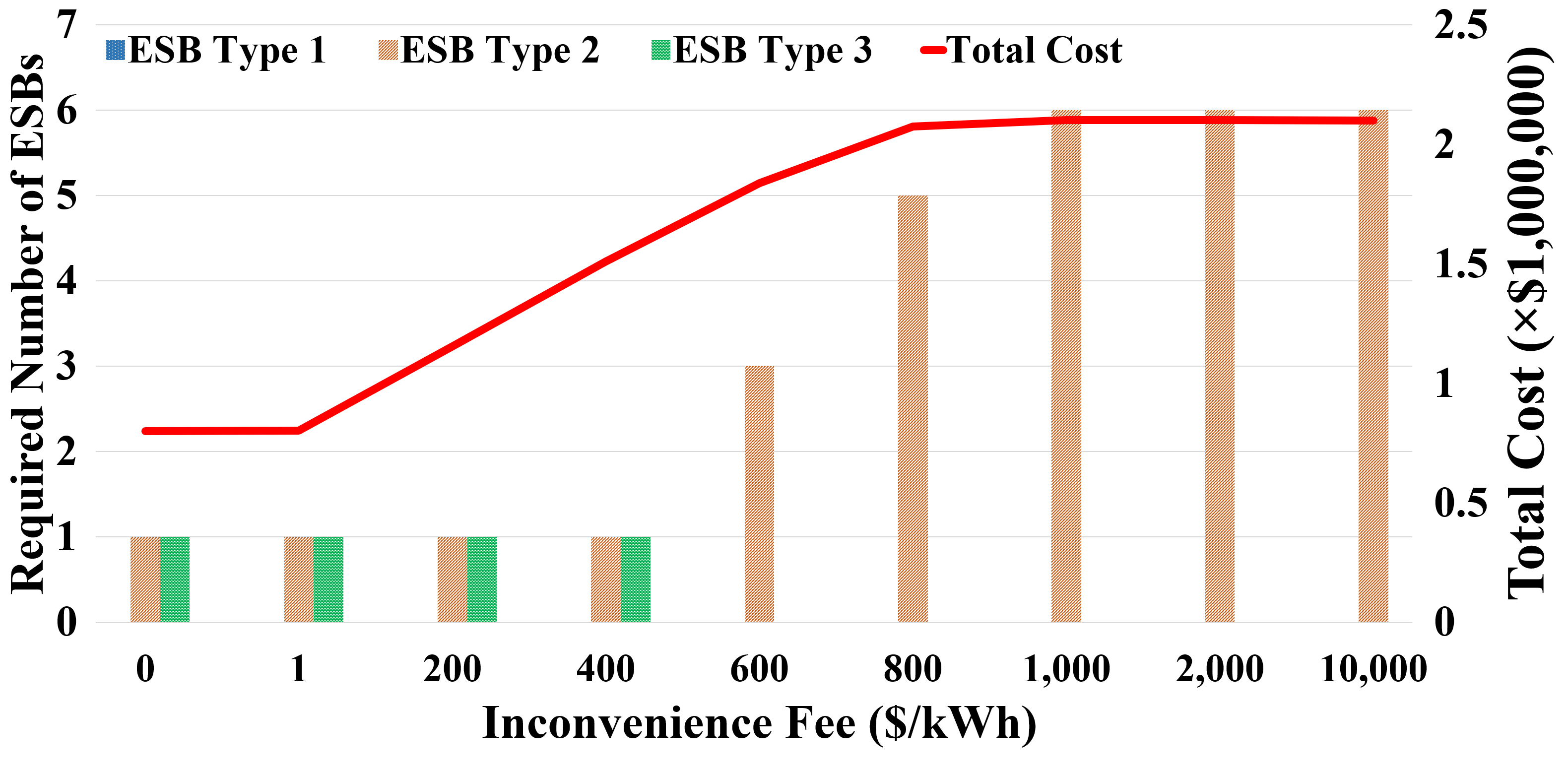}
        \caption{Problem instance 1-1-32.}
        \label{fig:Effect_inconv_on_Cost_ESB_2}
    \end{subfigure}
    \begin{subfigure}[b]{0.48\textwidth}
        \centering
        \includegraphics[scale=0.11]{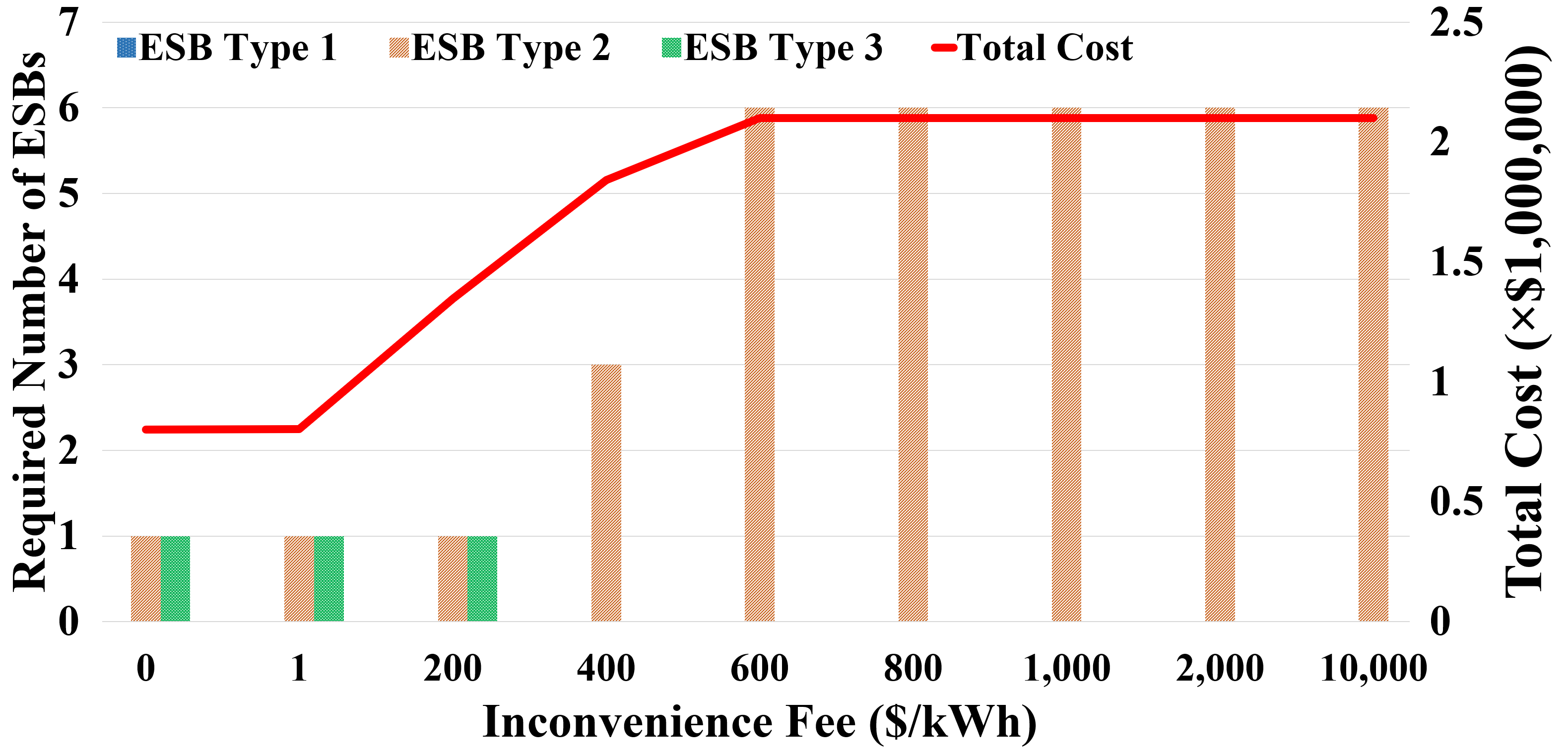}
        \caption{Problem instance 1-1-48.}
        \label{fig:Effect_inconv_on_Cost_ESB_3}
    \end{subfigure}
    \caption{Effect of inconvenience fee on the total cost and the required fleet composition of ESBs.}
    \label{fig:Effect_inconv_on_Cost_ESB_2_3}
\end{figure}

Figure \ref{fig:Effect_inconv_on_Cost_ESB_2_3} shows the effect of increasing the inconvenience fee of demand shifting on the total cost and the required fleet composition of the problem instances 1-1-32 and 1-1-48.

\end{document}